\documentclass [12pt]{article}
\usepackage{graphicx,amssymb,amsfonts,latexsym,amsmath,amsthm,times}
\usepackage{fancyhdr}
\usepackage{color}
\usepackage{microtype}
\usepackage[english]{babel}
\numberwithin{equation}{section}

\nonstopmode

\usepackage[utf8]{inputenc}

\usepackage{stmaryrd}

\usepackage{verbatim}
\usepackage{url}
\usepackage[numbers,sort]{natbib}
\usepackage{natbibspacing}
\usepackage{microtype}

\usepackage{subfigure}

\usepackage{tikz}
\usepackage{pgfmath}

\usetikzlibrary{positioning}
\usetikzlibrary{calc}
\usetikzlibrary{chains}
\usetikzlibrary{shapes.misc}
\usetikzlibrary{shapes.symbols}
\usetikzlibrary{shapes.geometric}
\usetikzlibrary{fit}
\usetikzlibrary{shadows}
\usetikzlibrary{arrows}
\usetikzlibrary{trees}
\usetikzlibrary{decorations.pathmorphing}

\pgfdeclarelayer{background}
\pgfdeclarelayer{foreground}
\pgfsetlayers{background,main,foreground}

\let\vref=\ref
\newcommand{\assign}{:=}

\newcommand{\jump}[1]{\left\llbracket#1\right\rrbracket}
\newcommand{\D}{\mathsf{D}}

\newcommand{\mathd}{\,\mathrm{d}}
\newcommand{\doi}[1]{\href{http://dx.doi.org/#1}{doi: #1}}

\usepackage{hyperref}

\begin{document}

\footnotesize {\flushleft \mbox{\bf \textit{Math. Model. Nat.
Phenom.}}}
 \\
\mbox{\textit{{\bf Vol. X, No. X, 2011, pp. 1-\pageref{pg:atend}}}}

\thispagestyle{plain}

\vspace*{2cm} \normalsize
\begin{center}
  \Large \bf Viscous Shock Capturing in a Time-Explicit Discontinuous~Galerkin~Method
\end{center}

\vspace*{1cm}

\centerline{\bf A. Klöckner$^a$\footnote{Corresponding author. E-mail: kloeckner@cims.nyu.edu},
  T.~Warburton $^b$ and J.~S.~Hesthaven$^c$ 
  }

\vspace*{0.5cm}

\centerline{$^a$ Courant Institute of Mathematical Sciences, New York University, New York, NY 10012}

\centerline{$^b$ Department of Computational and Applied Mathematics, Rice University, Houston, TX 77005 }

\centerline{$^c$ Division of Applied Mathematics, Brown University, Providence, RI 02912}


\vspace*{1cm}

\noindent {\bf Abstract.}
We present a novel, cell-local shock detector for use with
discontinuous Galerkin (DG) methods. The output of this detector is a
reliably scaled, element-wise smoothness estimate which is suited as a
control input to a shock capture mechanism. Using an artificial
viscosity in the latter role, we obtain a DG scheme for the numerical
solution of nonlinear systems of conservation laws.  Building on work
by Persson and Peraire, we thoroughly justify the detector's design
and analyze its performance on a number of benchmark problems. We
further explain the scaling and smoothing steps necessary to turn the
output of the detector into a local, artificial viscosity. We close by
providing an extensive array of numerical tests of the detector in
use.

\vspace*{0.5cm}




\label{cha:viscous}
\section{Introduction}
\label{sec:viscous-intro}

Discontinuous Galerkin methods
\cite{reed_triangular_1973,lesaint_finite_1974,cockburn_runge-kutta_1990,hesthaven_nodal_2007}
are a high-order accurate, geometrically flexible, and robust means of
approximating solutions of systems of hyperbolic conservation laws.
For linear conservation laws, these schemes easily deliver highly
accurate solutions without much effort. For nonlinear hyperbolic
systems, the situation is more complicated.  If the solution of the
system remains smooth for the entire time under consideration, and if
thereby the decay of modal coefficients is fast enough, the method may
be used with little modification for a so-called ``nodal approach''.
Optionally, aliasing error in the computation of integrals for
stiffness and mass matrices can be avoided by the introduction of
quadrature schemes of sufficient order \cite{hesthaven_nodal_2007}.

If however the solution does not stay smooth for long enough periods
of time, the arising discontinuities pose a number of problems which
have been the subject of intense study since the early days of
scientific computation and numerical analysis. \citep[e.g.][and
references therein]{gottlieb_gibbs_1997} Our goal is to seek out a
method that is able to reliably detect the occurrence of Gibbs
phenomena (which represent the main issue with the discontinuous
solution) in the context of the discontinuous Galerkin method. In this
paper, the subsequent mitigation of then phenomenon is then achieved
through a simple artificial viscosity.

Many authors have proposed methods to capture shocks within a DG
setting, by different methods.  \emph{Flux limiting}, which has been
both successful and popular with Finite Volume practitioners, was
combined with DG immediately in conjunction with the resurgence of
interest in the method in the late 1980s. \citep{cockburn_tvb_1989,
cockburn_tvbsys_1989, cockburn_runge-kutta_1990,
cockburn_runge-kutta_1998,burbeau_problem-independent_2001,
dolejsi_aspects_2003, tu_slope_2005, kuzmin_flux_2005,
krivodonova_limiters_2007, xu_hierarchical_2009}. A common theme to
limiting is that the solution is modified in some way to retain
desirable properties such as positivity and freedom from spurious
oscillation, and in doing so, reaches various (often low) orders of
accuracy.

\emph{Artificial viscosity} methods, on the other hand, take the position
that the only hope of resolving a discontinuity by a high-order
approximation lies in smoothing it out. These methods date back to
\citet{neumann_method_1950}, were first used in the context of finite
difference methods \citep{lapidus_detached_1967}, and then spread into
finite element literature (see, e.g., the study by
\citet{john_finite_2008} for a review) and were also applied to
time-dependent discontinuous Galerkin methods very early on
\citep{bassi_accurate_1994}, and have since enjoyed continuing
popularity \citep[e.g.][]{burman_nonlinear_2007}.

One obvious improvement on \emph{global} artificial viscosities is a
more selective application of smoothing, guided by a detector.  There
has been a recent resurgence of interest in such approaches
\cite{persson_subcell_2006,barter_shock_2009} in the context of DG.
The methods discussed in this article aim to improve on these latter
schemes, where we would like to emphasize that our detector is
\emph{not} intrinsically tied to guiding the application of an
artificial viscosity. With appropriate rescaling, it might be suitable
in a multitude of other scenarios requiring shock detection.

Other variants of artificial viscosity methods exist as well. The
method of \emph{Spectrally Vanishing Viscosity}
\citep[e.g.][]{tadmor_convergence_1989,kirby_stabilisation_2006} is
similar in spirit, but tries to restrict its smoothing action to
high-frequency solution components.

One final approach of dealing with discontinuities is that of adapting
the mesh and adding resolution.  It is generally thought that
`shocks', i.e. genuine discontinuities, do not exist in nature
\citep{woodward_numerical_1984}, and thereby, if only enough
resolution were available, the problem of shock capturing would vanish
by itself. While nature may obey this statement, mathematical models
of it often do not (e.g. Burgers' equation), and so one needs to
``help a little''--for example by adding an artificial viscosity
\cite[e.g.][]{hartmann_adaptive_2006}. Others contend that the wiggles
are worth keeping simply as indicators of numerical trouble
\citep{gresho_dont_1981}. Further, while adaptivity certainly is a
useful technique in capturing shocks
\cite{warburton_galerkin_1999,flaherty_adaptive_1997,xu_high-order_2010},
it too depends on a detector that reliably tells the method where
refinement is necessary.

\subsection{The Discontinuous Galerkin Method}

\label{sec:dg-overview}

\begin{figure}
\begin{center}
\begin{tikzpicture}[
  densemat/.style={fill=gray!30},
  ]
  \draw [densemat] (-1,-1) rectangle (1,1) ;
  \draw (2,-1) rectangle +(4,2) ;
  \draw [densemat] (-2,-1) rectangle (-5,1) ;

  \draw [densemat] (2,1) rectangle +(1,-1) ;
  \foreach \ys/\ye in {0.1/0.3, -0.25/-0.7, 0.6/0.8}
  { \draw [densemat] (3,\ys) rectangle (4,\ye) ; }
  \foreach \ys/\ye in {0.4/0.1, -0.1/-0.3, -0.9/-1, 1/0.9}
  { \draw [densemat] (4,\ys) rectangle (5,\ye) ; }
  \foreach \ys/\ye in {0.8/0.6, 0.3/0, -0.9/-0.7, -0.2/-0.5}
  { \draw [densemat] (5,\ys) rectangle (6,\ye) ; }

  \node at (2.5,0.5) {$M^{k,A_1}$} ;
  \node at (3.5,-0.475) {$M^{k,A_2}$} ;
  \node at (4.5,0.25) {$M^{k,A_3}$} ;
  \node at (5.5,-0.35) {$M^{k,A_4}$} ;

  \node at (0,0) {$(M^k)^{-1}$} ;
  \node at (1.5,0) {$\cdot$} ;
  \node at (-1.5,0) {$=$} ;
  \node at (-3.5,0) {$L^k$} ;

  \draw [|<->|] (6.25,-1) -- +(0,2) node [pos=0.5, anchor=west] {$N_p$} ;

  \draw [|<->|] (2,-1.25) -- (3,-1.25) node [anchor=west] {$N_{fp}$} ;
\end{tikzpicture}
\end{center}
\caption{Construction of the Lifting Matrix $L^k$.}
\label{fig:lifting-matrix}
\end{figure}
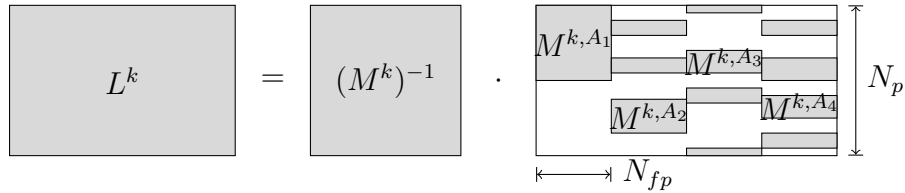

Discontinuous Galerkin (DG) methods
\cite{reed_triangular_1973,lesaint_finite_1974,cockburn_runge-kutta_1990,hesthaven_nodal_2007}
are a combination of ideas from Finite-Volume and Spectral Element methods.
We consider DG methods for the approximate solution of a hyperbolic system of
conservation laws
\begin{equation} 
  u_t + \nabla \cdot F (u) = 0 \label{eq:claw} 
\end{equation}
on a domain $\Omega = \biguplus_{k=1}^K \D_k \subset \mathbb R^d$ consisting of
disjoint, face-conforming tetrahedra $\D_k$ with boundary conditions
\[ 
  u|_{\Gamma_i} (x, t) = g_i (u (x, t), x, t), \hspace{2em} i = 1, \ldots, b,
\]
at inflow boundaries $\biguplus \Gamma_i \subseteq \partial \Omega$. 
We find a weak form of (\ref{eq:claw}) on each element $\D_k$:
\[
  0 = \int_{\D_k} u_t \varphi + [\nabla \cdot F (u)] \varphi \mathd x
   = \int_{\D_k} u_t \varphi - F (u) \cdot \nabla \varphi \mathd x
  + \int_{\partial \D_k} ( \hat{n} \cdot F)^{\ast} \varphi \mathd S_x, 
\]
where $\varphi$ is a test function, and $( \hat{n} \cdot F)^{\ast}$ is a
suitably chosen numerical flux in the unit normal direction $\hat{n}$.
Following \cite{hesthaven_nodal_2007}, we may find a so-called `strong'-DG form of
this system as
\begin{equation}
  0 = \int_{\D_k} u_t \varphi + [\nabla \cdot F (u)] \varphi \mathd x-
  \int_{\partial \D_k} [ \hat{n} \cdot F - ( \hat{n} \cdot F)^{\ast}] 
  \varphi \mathd S_x.
  \label{eq:strong-dg}
\end{equation}
by integrating by parts once more.  We seek to find a numerical vector
solution $u^k \assign u_N|_{\D_k}$ from the space $P_N^n (\D_k)$ of
local polynomials of maximum total degree $N$ on each element, where
$n$ is the number of equations in the hyperbolic system
\eqref{eq:claw}.  We choose the scalar test function $\varphi \in P_N
(\D_k)$ from the same space and represent both by expansion in a basis
of $N_p\assign \operatorname{dim} P_N(\D_k)$ Lagrange polynomials
$l_i$ with respect to a set of interpolation nodes
\cite{warburton_explicit_2006}. We define the mass, stiffness,
differentiation, and face mass matrices
\begin{subequations} 
  \label{eq:dg-global-matrices} 
  \begin{align} 
    M_{i j}^k & \assign \int_{\D_k} l_i l_j \mathd x, &
    S_{i j}^{k, \partial\nu} & \assign \int_{\D_k} l_i \partial_{x_{\nu}} l_j \mathd x,\\ 
    D^{k, \partial\nu} & \assign (M^k)^{- 1} S^{k, \partial\nu},&
    M_{i j}^{k, A} & \assign \int_{A \subset \partial \D_k} l_i l_j \mathd S_x.
  \end{align} 
\end{subequations}
Using these matrices, we rewrite \eqref{eq:strong-dg} as
\begin{align} 
  0 & =  M^k \partial_t u^k 
  + \sum_{\nu} S^{k, \partial_{\nu}} [F(u^k)]
  - \sum_{F \subset \partial \D_k} M^{k, A} [
  \hat{n} \cdot F - ( \hat{n} \cdot F)^{\ast}], \notag \\
  \label{eq:semidiscrete-dg} 
  \partial_t u^k & = 
  - \sum_{\nu} D^{k, \partial_{\nu}} [F(u^k)]
  + L^k [ \hat{n} \cdot F 
  - ( \hat{n} \cdot F)^{\ast}] |_{A \subset \partial \D_k}.
\end{align}
The matrix $L^k$ used in \eqref{eq:semidiscrete-dg} deserves a little
more explanation. It acts on vectors of the shape
$[u^k|_{A_1},\dots,u^k|_{A_4}]^T$, where $u^k|_{A_i}$ is the vector of
facial degrees of freedom on face $i$. For these vectors, $L^k$
combines the effect of applying each face's mass matrix, embedding the
resulting facial values back into a volume vector, and applying the
inverse volume mass matrix.  Since it ``lifts'' facial contributions
to volume contributions, it is called the \emph{lifting matrix}. Its
construction is shown in Figure \ref{fig:lifting-matrix}. 

It deserves explicit mention at this point that the left
multiplication by the inverse of the mass matrix that yields the
explicit semidiscrete scheme \eqref{eq:semidiscrete-dg} is an
element-wise operation and therefore feasible without global
communication. This strongly distinguishes DG from other finite
element methods. It enables the use of explicit (e.g., Runge-Kutta)
time stepping and greatly simplifies parallel implementation efforts.

\section{Basic Design Considerations}
\label{sec:viscous-design-considerations}

This article describes a method for detecting (and also capturing)
shocks in the context of DG methods. One particular motivation for us
was our recent work on the efficient mapping of DG onto massively
parallel throughput-oriented computer architectures
\citep{kloeckner_2009}, where we demonstrated a method to quickly
compute the vector $A(x)$ for a linear discontinuous Galerkin operator
$A$ and a state vector $x$ using graphics hardware (i.e. Graphics
Processing Units or ``GPUs''). The present article describes one stepping stone
on the way to generalizing the applicability of GPU-DG to nonlinear
problems.

In briefly explaining the unique environment present on GPUs, we seek
to inform the reader on the considerations that guided our approach.
On wide-SIMD, parallel architectures such as the GPUs of
\citep{kloeckner_2009}, memory is at a premium and scattered memory
access is particularly expensive. As a consequence, we argue
that \emph{matrix-free} methods such as the one of
\citep{kloeckner_2009}, if they can be implemented efficiently, will
always hold a significant performance advantage over approaches that
have to build, keep in memory, and constantly access a pre-built
sparse matrix for $A$, because such a computation is necessarily bound
by the speed at which matrix entries can be streamed into the core,
where they are then used exactly once and discarded
\citep{bell_spmv_2008}. A matrix-free approach has far more freedom to
exploit local structure and re-use data. We will therefore focus our
investigation on matrix-free methods.

This choice has important ramifications. One consequence of it affects
the trade-off by which one chooses between implicit and explicit time
stepping. Consider the case of implicit time integrators, in which one
must constantly solve large linear systems of equations. Direct,
factoring solvers for sparse matrices are as yet unavailable on
massively parallel hardware, and even if they were, they would doubly
suffer from the issues that sparse matrices encounter. One therefore
naturally looks towards iterative methods for solving large sparse
systems. For the complicated linearized systems arising from the
nonlinear hyperbolic conservation laws we are targeting in this
article, these methods generally need help in the form of a
preconditioner in order to be efficient. This is the next implication
of the choice of matrix-free methods: One automatically chooses to not
use the substantial body of literature showing how a preconditioner may
be built from a known sparse matrix. Instead, one needs to invest
further work designing and testing preconditioners (using e.g.
multi-grid or domain-decomposition methods), and, in addition to the
design time spent, these preconditioners may carry significant
additional computational expense, typically through their communication
needs.  In addition, Krylov methods (which are frequently used to
solved the arising large, sparse linear systems) in particular
involve global reductions (in the form of inner products) which are
known to not achieve peak performance on graphics processors
\citep{harris_optimizing_2007}.  Worse, the nonlinear PDEs we are
targeting in this paper require a nonlinear system of equations to be
solved (likely by Newton iteration, which in turn requires Jacobians to
be evaluated).

This collection of drawbacks and uncertainties in the application of
implicit time integration on massively parallel hardware makes it seem
opportune to examine the use of explicit time steppers, which were
already used with good success in \citep{kloeckner_2009}. We aim to
find out if the single big disadvantage of explicit methods, namely
their small time step restriction, can be offset by the judicious
choice of methods combined with the advantages conferred by the
hardware.

Since the scheme we are aiming to design involves the use of artificial
viscosity, the scaling of the explicit time step is typically given by
\begin{equation}
  \Delta t \sim \frac 1 {\displaystyle \lambda_{\text{max}} \frac{N^2}{h} +
  \|\nu\|_{L^\infty} \frac{N^4}{h^2}},
  \label{eq:viscous-timestep-restriction}
\end{equation}
where $\lambda_{\text{max}}$ is the largest characteristic velocity
and $\nu$ is the magnitude of the viscosity, $h$ is the local mesh
size and $N$ is the approximation's polynomial degree
\cite{hesthaven_nodal_2007}.
Within \eqref{eq:viscous-timestep-restriction}, the numerical
diffusion time scale $N^4\|\nu\|_{L^\infty} / h^2$ can be rather
damaging, as it contains discretization-dependent factors at high
exponents.

Luckily, \eqref{eq:viscous-timestep-restriction} does not tell the
entire story. For example, we expect the occurrences of high
viscosity $\nu$ to be localized in both space and time. 
Localization in space could conceivably be dealt with using local time
stepping, but this is beyond the scope of this article. Localization in time
on the other hand is easily dealt with by the use of
time-adaptivity \citep[e.g.][]{dormand_family_1980}.  Adaptivity in
time is particularly important for explicit time stepping of
artificial-viscosity-enhanced PDE solvers.

One further aspect of the time discretization should be considered: Much of the
effort in this research is targeted at mitigating the effect of oscillations in
the spatial discretization of a conservation law that trace their roots back to
the polynomial expansions used for them. Time discretizations, however, are
equally based on polynomials, and many varieties of so-called Strong Stability
Preserving (SSP) time integrators have been devised to mitigate oscillations
originating in temporal expansions \citep{shu_tvd_1988}.  Even embedded pairs
of SSP Runge-Kutta methods are available \citep{gottlieb_strong_2011}. Based on
initial experiments, it appears that in the setting of this work,
spatially-generated oscillations by far dominate their temporal cousins at the
time step sizes encountered. Thus varying the time integration method does not
have an appreciable effect on the reported results.

In summary, the emergence of massively parallel hardware along with the
use of purposefully chosen, adaptive time discretizations may help
explicit methods be competitive with implicit methods for the
integration of large-scale nonlinear systems, a few of which we will
introduce next.

\section{Applications and Equations}
\label{sec:viscous-equations}
\subsection{Advection Equation}
\label{sec:viscous-eqn-advection}
At the very simple end of the spectrum of hyperbolic conservation laws, the
\emph{advection equation}
$\partial_t u+v \cdot \nabla_x u = 0$
transports its initial condition along its one characteristic,
described by the velocity vector $v$. We will apply artificial
viscosity to this PDE as
\[\partial_t u+ v\cdot \nabla_x u = \nabla_x\cdot (\nu \nabla_x u).\]
Here, and in all further equations, it is important to write the
viscosity in ``conservation'' form $\nabla_x\cdot (\nu \nabla_x u)$.
The desired consequence of this is that the resulting DG method will
be conservative \cite{arnold_unified_2002}.

In DG discretizations of this equation, we use a conventional upwind flux
in a strong-form DG formulation.
The diffusion term $\nabla_x\cdot (\nu \nabla_x u)$ is discretized by
a first-order (``dual'') \emph{interior penalty method}
\citep{arnold_unified_2002}, with the gradient being computed in
strong form, and the divergence computed in weak form. The
diffusive fluxes are given by
\[
  u_N^\ast:= \{u_N\},\qquad
  \sigma_N^\ast:= \{\nu \nabla_{x,h} u_N\}-\frac{N^2}h \nu \jump{ u_h},
\]
where $\sigma_N$ is the discretization of $\nu\nabla_x u$.
\subsection{Second-Order Wave Equation}
\label{sec:viscous-eqn-wave}
The wave equation $\partial t^2 u+c^2 \triangle u=0$ is valuable for testing
artificial viscosity methods because it is the simplest system where the
effects of two coupled characteristics may be observed.  We rewrite this PDE as
a first-order system of conservation laws and apply artificial viscosity to
this system to obtain
\begin{subequations}
  \label{eq:viscous-wave-eqn}
  \begin{align}
  \partial_t u+c \nabla_x \cdot v &= \nabla_x\cdot(\nu \nabla_x u),\\
  \partial_t v+c \nabla_x u &= \nabla_x\cdot (\nu \nabla_x v),
  \end{align}
\end{subequations}
where we have again been careful to use the conservative form of the
diffusive term. The vector diffusion term $\nabla_x\cdot (\nu \nabla_x
v)$ is to be read as the diffusion $\nu$ being applied to each
component separately.  The discontinuity sensor to be described below
operates on the scalar component $u$. In DG discretizations of this
equation, we again use a conventional upwind flux in a strong-form DG
formulation.  The diffusion terms are discretized in analogy to the
preceding section.

\subsection{Euler's Equations of Gas Dynamics}
\label{sec:viscous-eqn-euler}
Lastly, the system of conservation laws that justifies the effort
spent on this study, \emph{Euler's equations of gas dynamics},
broadly applies to compressible, inviscid flow problems.
As in Section \ref{sec:viscous-eqn-wave}, we are again choosing to use a
single artificial viscosity $\nu$ that applies to all components, such 
that we get the viscosity-endowed system
\begin{subequations}
\begin{align}
\partial_t\rho+ \nabla_x\cdot(\rho\mathbf u)
  &=\nabla_x \cdot (\nu \nabla_x \rho),\\
\partial_t(\rho\mathbf{u})+ \nabla_x\cdot(\mathbf u\otimes(\rho \mathbf u))+\nabla_x p
  &=\nabla_x \cdot (\nu \nabla_x (\rho \mathbf u)),\\
\partial_t E+ \nabla_x\cdot(\mathbf u(E+p))
  &=\nabla_x \cdot (\nu \nabla_x E).
\end{align}
\end{subequations}
The discontinuity sensor to be described below operates on the component
$\rho$.  In contrast to \citet{persson_subcell_2006}, we find that a
Navier-Stokes-like physical viscosity provides insufficient control of
oscillations in $\rho$. 

In DG discretizations of this system, a \emph{local
Lax-Friedrichs} (or \emph{Rusanov}) flux 
\[
  \hat n \cdot F_N^{\ast}\assign
  \hat n\cdot \frac{F(u^+)+F(u^-)}2
  - \frac{\lambda_{\text{max}} }2(u^+-u^-),
\]
in weak-form DG is commonly used. The diffusion term is discretized as in
Section \ref{sec:viscous-eqn-wave}.  As above, a quadrature exact to degree
$3N$ is used to integrate the nonlinearity.

\section{A Smoothness-Estimating Detector}
\label{sec:viscous-construction}

Detectors for the selective application of artificial viscosity have
been built in a large variety of ways. The most popular, perhaps, is
sensing on the $L^2$ norm of the residual of the variational form
\citep{bassi_accurate_1994,jaffre_convergence_1995}.
\citet{hartmann_adaptive_2006} employs a similar indicator that
includes sensing of the primary orientation of the discontinuity and
performs anisotropic mesh refinement based on this data.

Other detectors in the literature employ information gathered not on
the whole volume of the domain, but only on element faces
\citep{bassi_high-order_1997}. Specializing further, some methods use
the magnitude of the facial inter-element jumps as an indicator of how
well-resolved the solution is and to what degree it has converged
\citep{barter_shock_2009,feistauer_robust_2007}.
A further approach to shock detection repurposes entropy pairs,
objects from the solution theory for scalar conservation laws, for the
purposes of shock detection \citep{guermond_entropy_2008}.

Our approach most directly traces its lineage to work by
\citet{persson_subcell_2006}, which addresses one crucial shortcoming
in much of the above work: scaling. Many of the quantities discussed
clearly relate directly to how well-resolved (and smooth) the
approximate solution of the system is. It is however rarely clear how
large a value of the quantity in question indicates that a problem
exists, and a variety of ad-hoc scaling choices are proposed, often by
the maximum of the quantity found across the domain, or by the
element-local norm, but without assigning an explicit meaning to the
scaled quantity.

The method by \citet{persson_subcell_2006} also performs scaling by
the element-local $L^2$ norm $\|q_N\|_{L^2(\D_k)}$ of the discretized
value of the quantiy $q_N$ to be sensed on. On each element $\D_k$, it
obtains a value
\begin{equation}
  \label{eq:viscous-pp-indicator}
  S_k:=\frac
  {(q_N,\phi_{N_p-1})^2_{L^2(\D_k)}}
  {\|q_N\|^2_{L^2(\D_k)}},
\end{equation}
where $\{\phi_n\}_{n=0}^{N_p-1}$ is an orthonormal basis for the
expansion space \citep[see e.g.][]{dubiner_spectral_1991,
koornwinder_two-variable_1975} numbered from 0. Simply put, $S_k$
reflects the (squared) fraction of $q_N$'s mass contained in the
highest mode of the expansion, relative to all mass present on the
element. \citet{persson_subcell_2006} then invoke an analogy to
Fourier expansions, where a continuous function (roughly) can be
recognized by having Fourier expansions in which the $n$th mode's
magnitude scales at most as $1/n^2$. In doing so, they have
conveniently solved the issue of scaling--it is now understood what
$S_k$ measures and what value it is supposed to take on for which
degree of smoothness.  Based on this analogy, they argue that $S_k$
should have a magnitude of $1/N^4$ for $q_N$ to be continuous, or,
alternatively, that smoothing by artificial viscosity should activate
if $S_n> 1/N^4$.

They achieve this activation through a sequence of mapping steps.
First, they take the logarithm
\[
  s_k\assign \log_{10} S_k
\]
to obtain a quantity that scales linearly with the decay exponent,
which they put in relation to a quantity $s_0$ that they claim should
scale as $1/N^4$. We believe this is a typographical error in their
paper, because for proper comparability, $s_0$ should scale with the
\emph{logarithm} of $1/N^4$. Through the application of a mapping,
they obtain the final per-element viscosity
\begin{equation}
  \label{eq:viscous-pp-activation-map}
  \nu_k(s_k) = \nu_0 \begin{cases}
    0 & s_k< s_0 - \kappa, \\

    \frac12\left(1+\sin\frac{\pi (s_k-s_0)}{2\kappa} \right)
      & s_0-\kappa \le s_k \le  s_0 + \kappa, \\

    1& s_0 + \kappa < s_k,
  \end{cases}
\end{equation}
where $\nu_0$ is the maximum viscosity, which \citet{persson_subcell_2006}
suggest to scale with $h/N$ and $\kappa$ is the width of the activation
``ramp''.

The focus of the remainder of this article is to identify a number of
issues and make a number of improvements to this method of finding an
artificial viscosity. 

\subsection{Estimating Solution Smoothness}
\label{sec:viscous-estimating-smoothness}
\def\estimationlegend{
  Subfigure (a) shows the nodal data and its unique polynomial
  interpolant. Subfigure (b) shows the modal coefficients of
  a Legendre expansion of the function in (a), the processing
  of these coefficients, and the unprocessed and postprocessed
  smoothness estimates.
}
\newcommand{\viscousshowestimation}[2]{
  \begin{figure}
    \begin{center}
    \subfigure[]{
      \label{fig:viscous-estimation-#1-spatial}
      \includegraphics[width=0.4\textwidth]{modeproc-#1-spatial.pdf}
    }
    \subfigure[]{
      \label{fig:viscous-estimation-#1-modal}
      \includegraphics[width=0.4\textwidth]{modeproc-#1-modes-pre.pdf}
    }
    \end{center}%
    \caption[#2]{#2
       \estimationlegend
    }%
    \global\def\estimationlegend{}%
    \label{fig:viscous-estimation-#1}%
  \end{figure}
}

Before we begin our discussion of the refinements to the method, let
us set the stage by discussing the type of numerical method at which
the to-be-designed detector is aimed. As was already discussed, for
methods of low approximation order (and polynomial degrees
$N\lessapprox 2$), the flux limiting literature provides plenty of
alternatives for shock capturing, and therefore will not be the main
target area for our work. Since our method, like the work of
\citet{persson_subcell_2006}, will try to extract smoothness
information from the modal expansion of the solution, it is our hope
that the expansion at these degrees already contains enough smoothness
information to be viable as a basis for an artificial viscosity.
Lastly, at degrees $N\gtrapprox 5$, there is guaranteed to be
sufficient smoothness information, though the time step restriction
\eqref{eq:viscous-timestep-restriction} may make these approximations
somewhat impractical.

We begin our deconstruction and rebuild of the Peraire-Persson estimator
by examining the assumption that, like for Fourier series, smoothness
can be estimated by modal decay. In Fourier series, this can be
justified by viewing what happens if a derivative of an expanded
function is taken (and hence smoothness is reduced)--the $n$th
coefficient's magnitude gets multiplied by $n$. This results in the
identity
\begin{equation}
  \label{eq:viscous-derivative-fourier}
  \left \| \frac {d}{dx} e^{inx} \right \|_{L^p((-\pi,\pi))} 
  = n \left\|e^{inx}\right\|_{L^p((-\pi,\pi))}
  \qquad \text{for $p\in[1,\infty]$}.
\end{equation}
A polynomial analog for \eqref{eq:viscous-derivative-fourier} 
is provided by Bernstein's inequality
\citep{warburton_taming_2008,borwein_polynomials_1995}
\begin{equation}
  \label{eq:viscous-bernsteins-ineq}
  \left | \frac {d}{dx} P(x) \right | 
  \le \frac n{\sqrt{1-x^2}} \|P(x)\|_{L^\infty([-1,1])}
  \qquad \text{for $P\in P^n([-1,1])$, $x\in[-1,1]$}.
\end{equation}
While it conveniently exhibits the same scaling as its Fourier
counterpart, unfortunately, this estimate breaks down near the domain
boundaries.  Markov's inequality \citetext{ibid.}
\begin{equation}
  \label{eq:viscous-markovs-ineq}
  \left \| \frac {d}{dx} P(x) \right \|_{L^\infty([-1,1])} 
  \le n^2 \left\|P(x)\right\|_{L^\infty([-1,1])}
  \qquad \text{for $P\in P^n([-1,1])$}.
\end{equation}
extends the estimate out to the domain boundary, at the expense of a
larger scaling. Further, it may be argued that if one wants to
transfer the knowledge gained from \eqref{eq:viscous-markovs-ineq} to
a modal setting, $L^\infty$ is the wrong norm, and one should consider
the $L^2$ norm instead to be able to benefit from Parseval's identity.
Fortunately, an $L^2$ analog of \eqref{eq:viscous-markovs-ineq} is
available \citep[and references therein]{warburton_taming_2008}
\begin{equation}
  \label{eq:viscous-l2-markovs-ineq}
  \left \| \frac {d}{dx} P(x) \right \|_{L^2([-1,1])} 
  \le \sqrt3 n^2 \left\|P(x)\right\|_{L^2([-1,1])}
  \qquad \text{for $P\in P^n([-1,1])$},
\end{equation}
known as an \emph{inverse inequality}. Taking into account
\eqref{eq:viscous-bernsteins-ineq} and
\eqref{eq:viscous-l2-markovs-ineq}, the polynomial analogy to the
Fourier case is therefore expected to carry over well for
non-smoothness occurring on the interior of each finite element,
whereas for non-smoothness at the domain boundary, the smoothness
measure will likely differ.

Having examined the viability of modal decay as an estimator for
smoothness, we seek to make the notion of modal decay more precise than
\eqref{eq:viscous-pp-indicator}. We presume that, for the modal
coefficients $\{\hat q_n\}_{n=0}^{N_p-1}$ of a member $q_N$ of the
$L^2$-orthonormal approximation space spanned by
$\{\phi_n\}_{n=0}^{N_p-1}$, modal decay is approximately representable
as
\begin{equation}
  \label{eq:viscous-modal-decay-model}
  |\hat q_n|\sim c n^{-s}.
\end{equation}
Taking the logarithm of the relationship \eqref{eq:viscous-modal-decay-model}
yields
\[
\log |\hat q_n|\sim \log(c) -s \log(n),
\]
an affine relationship whose coefficients $s$ and $log(c)$ may be
found through least-squares fitting, satisfying
\begin{equation}
  \label{eq:viscous-mode-least-squares}
  \sum_{n=1}^{N_p-1} \left|\log |\hat q_n|-(\log(c) -s \log(n))\right|^2
  \rightarrow \min!
\end{equation}
Observe that the decay rate of \eqref{eq:viscous-modal-decay-model}
has rather little to do with the presumed magnitude of the remainder
term of an expansion, on which most a-priori error estimates for
finite element solutions are based--these \emph{start} with an
assumption of sufficient smoothness.  There is a connection, however.
\citet{mavriplis_adaptive_1994}, in the context of mesh adaptation,
has used a similar least-squares fit to the modal decay, defining a
continuous function $\hat q(n)$ through the found fit.  She then
proceeds to estimate the remainder term of the expansion as
\[
  \| q- q_N\|_{L^2(\D^k)}^2 \approx
  \left(
  \frac{\hat q_N^2}{\frac{2N+1}{2}}
  +\int_{N+1}^\infty \frac{\hat q(n)^2}{\frac{2n+1}{2}}\mathd n
  \right).
\]
In a similar vein, Houston and Süli in \cite{houston_note_2005} use an $l^2$
fit like \eqref{eq:viscous-mode-least-squares} as a criterion for $hp$-adaptive
refinement. They obtain the approach from a discussion of results in
approximation theory \cite{davis_interpolation_1963}.

\viscousshowestimation{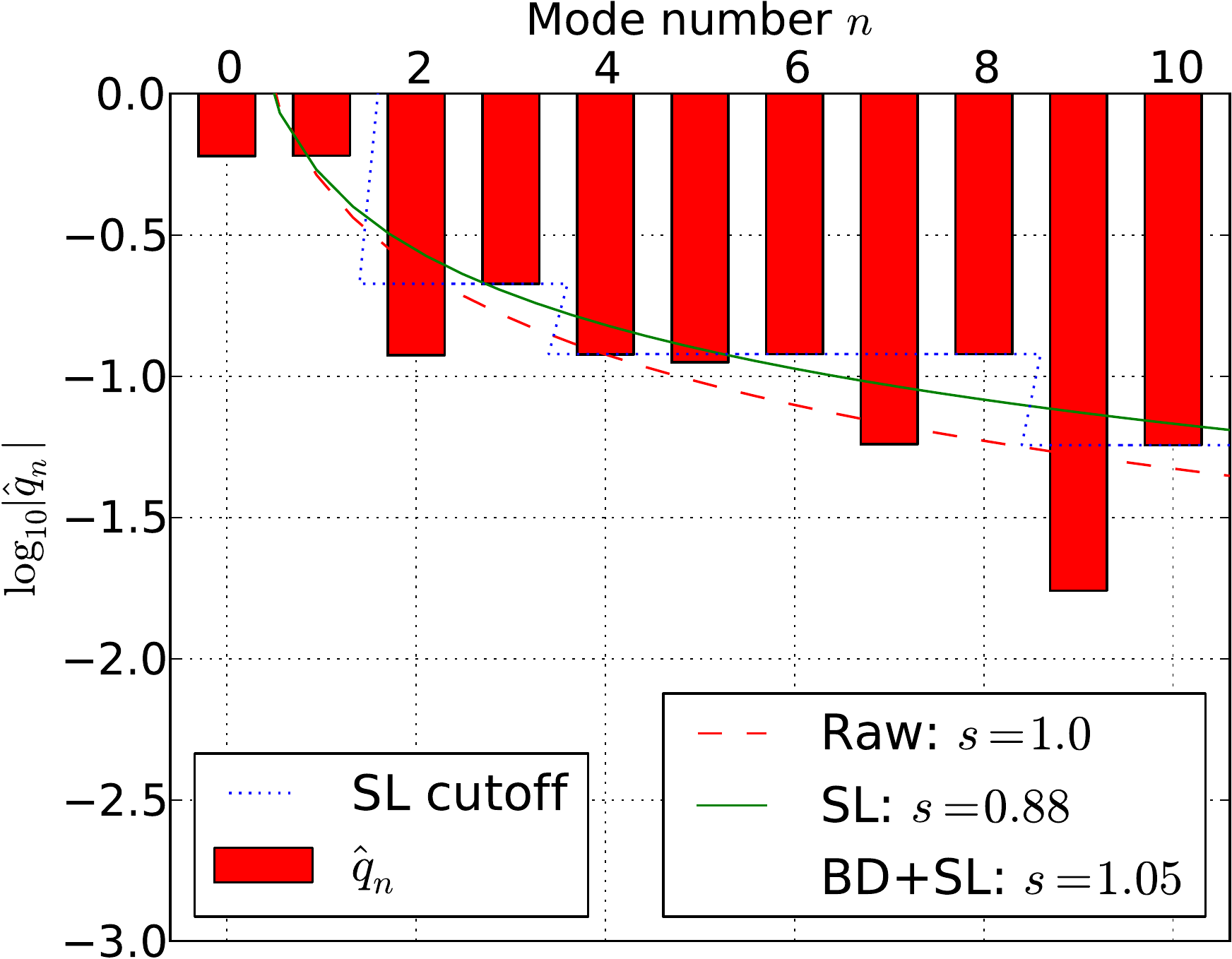}{
  Modal portrait for an approximant of a (discontinuous) 
  Heaviside jump function.
}

The least-squares procedure \eqref{eq:viscous-mode-least-squares} yields an
estimate $s$ of the decay exponent. If the analogy with Fourier modal decay
holds up, one would then expect $s\approx 1$ for a discontinuous $q$,
$s\approx2 $ for $q\in C^0\setminus C^1$, $s\approx 3$ for $q\in C^1\setminus
C^2$, and so forth. Figure \vref{fig:viscous-estimation-jump} shows a first
attempt at determining whether this is really the case by examining an
interpolant of a Heaviside jump function $H$ as shown in Figure
\ref{fig:viscous-estimation-jump-spatial}.  Figure
\ref{fig:viscous-estimation-jump-modal} shows the magnitudes of the first ten
modal coefficients along with the fitted curve (the dashed red line). The
obtained decay exponent $s$, shown in the legend next to the dashed red line,
matches the expectation well, giving a value of exactly 1.

\viscousshowestimation{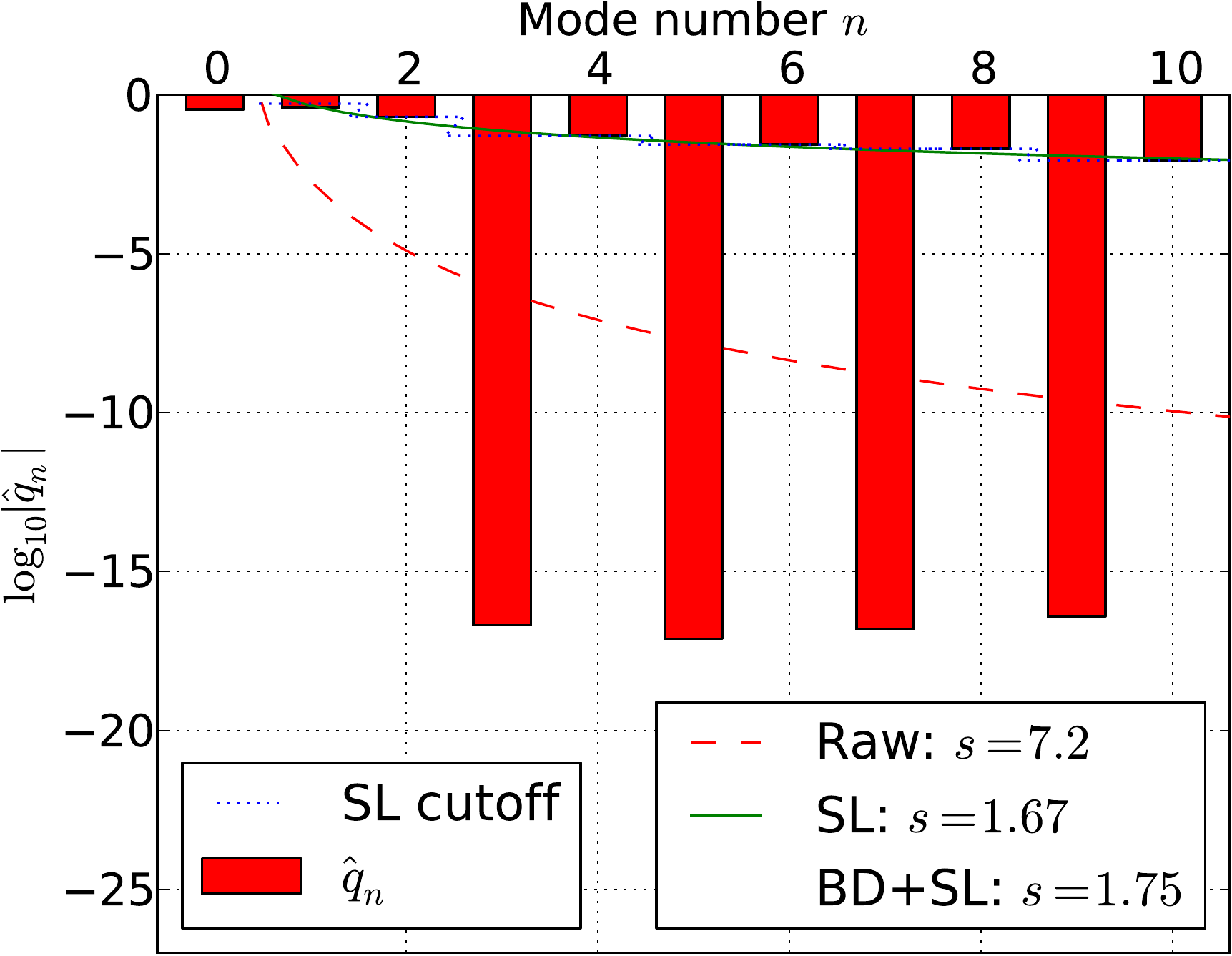}{
  Modal portrait for an approximant of a $C^0$ non-differentiable
  ``kink'' function.
}

Continuing this line of experimentation, we would like to move on to an
interpolant of a ``kink'' function $q(x):=xH(x)$.
The same observations as for the Heaviside function are shown in Figure
\vref{fig:viscous-estimation-kink}. Unfortunately, the figure reveals a
rather powerful shortcoming of the modal fit method as developed so
far. An odd-even effect draws the coefficients for the odd modes of
number three and greater to zero, leading to machine zeros ($\approx
10^{-15}$) in those approximate coefficient numbers. These ``fully
converged'' coefficients fool the estimator into an anomalous estimate of
far more smoothness than is actually present, leading to an
estimated decay exponent of about seven--far too high.


It is unfortunate that the fit can be misled that easily, but a close
look at Figure \ref{fig:viscous-estimation-kink-modal} will have
already revealed to the attentive reader that this is an easily
recoverable issue. Realize that the fit tries to model modal decay,
i.e. the shrinking of modal coefficient magnitudes $|\hat q_n|$ as $n$
increases. The model \eqref{eq:viscous-modal-decay-model} that is
fitted to the decay only generates monotone modal decays.  Figure
\ref{fig:viscous-estimation-kink-modal} is characterized by a strongly
non-monotone mode profile, and this is precisely what is misleading
the estimator. Consider this: Given a mode $n$ with a small
coefficient $|\hat q_n|$, if there exists another coefficient with
$m>n$ and $|\hat q_m|\gg |\hat q_n|$, then the small coefficient $|\hat
q_n|$ was likely spurious, just like the near-zero coefficients in
Figure \ref{fig:viscous-estimation-kink-modal} were spurious.
These spurious coefficients should hence be eliminated from the fit,
and this is what a new procedure, termed \emph{skyline pessimization},
achieves. From the modal coefficient magnitudes $\{|\hat
q_n|\}_{n=0}^{N_p-1}$, it generates a new set of modal coefficients by
\begin{equation}
  \label{eq:viscous-skyline}
  \bar q_n\assign
  \max_{i\in\{\min(n,N_p-2),\dots, N_p-1\}} |\hat q_i|
  \qquad \text{for $n\in\{1,2,\dots, N_p-1\}$}.
\end{equation}
The effect of the procedure is that each modal coefficient is raised
up to the largest higher-numbered modal coefficient, eliminating
non-monotone decay. Since odd-even effects in modal portraits (such as
the one of Figure \ref{fig:viscous-estimation-kink-modal} are a common
phenomenon, there is a slight modification in
\eqref{eq:viscous-skyline} accounting for the last mode, which is
forced to also be larger than the second-to-last mode. This would
become an issue if, for example, only the first nine modes of Figure
\ref{fig:viscous-estimation-kink-modal} were used, in which case the
smallness of the last coefficient would again cause an artificially
high smoothness exponent. Once skyline pessimization has been
performed, decay estimation \eqref{eq:viscous-mode-least-squares} is
applied to them in the same fashion as above, yielding a corrected
decay estimate.

The effect of skyline pessimization is shown in the modal portrait of
Figure \ref{fig:viscous-estimation-kink-modal} as a zig-zagged blue
line that appears to ``truncate'' the bars representing modal
coefficients at the level of the largest higher-numbered coefficient.
Further, the fitted decay curve is shown in green, along with the
resulting estimated decay exponent, labeled as ``SL''.  With skyline
pessimization in place, the estimated smoothness exponent for the
``kink'' example becomes 1.67--reasonably close to the expected value
of 2.

\viscousshowestimation{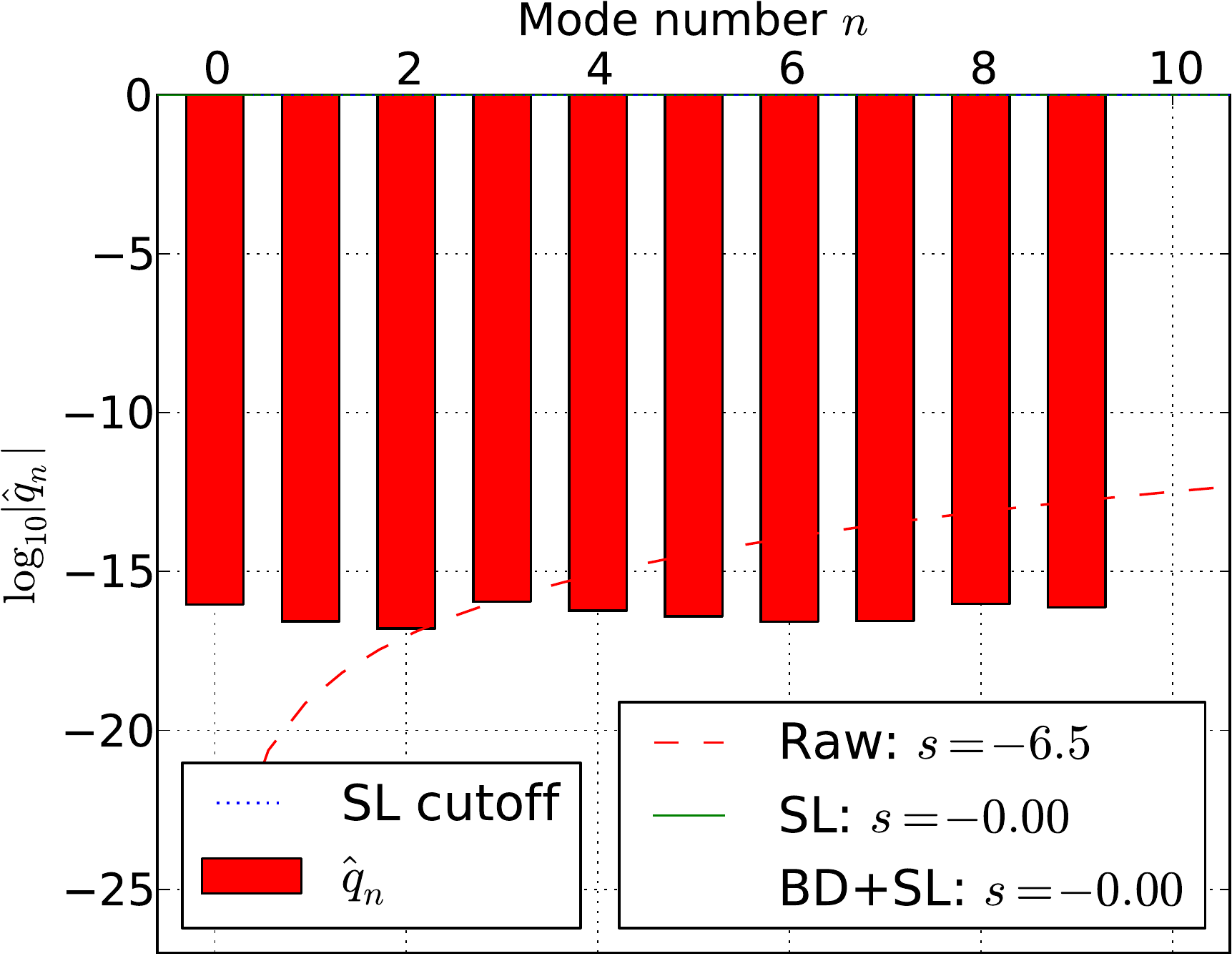}{
  Modal portrait for a function consisting of only the highest
  representable Legendre mode $\phi_{N_p-1}$ in an expansion of
  length 10.
}

The next-smoothest test of the estimator we consider is a truncated polynomial
$q(x)\assign x^2 H(x)$.
Obviously, $q\in C^1\setminus C^2$. As in the ``kink'' case, the modal decay
exhibits a pronounced odd-even discrepancy (not shown), that leads to
spuriously high ``raw'' smoothness exponent estimate of about 13. After skyline
pessimization, the estimate assumes nearly exactly the expected value, three.
The three artificial tests conducted so far confirm the premise on which the
estimator is built, namely that the smoothness of a function represented by a
Legendre expansion can be accurately estimated solely by examining its
coefficients.


By presenting a number of further tests, we hope to clarify the
behavior of the estimator as designed so far. A particularly
interesting case is shown in Figure
\vref{fig:viscous-estimation-last-mode}, which shows the estimator
applied to the highest mode present in the Legendre expansions of
length 10 which we have been considering. In a sense, this is the most
oscillatory, and thereby the least smooth, function that the expansion
can express. After skyline pessimization, this function is assigned a
smoothness exponent of zero--which in a Fourier setting would correspond
to white noise.


The next two tests are concerned with very smooth functions
($\cos(3+\sin(1.3x))$ and $\sin(\pi x)$) and confirm that the estimator
recognizes them as such.  While the smoothness values (both around four)
assigned to them are not as meaningful as the results in the low-smoothness
examples, this is not necessarily a problem. As long as the estimator can
sharply pick up non-smoothness on a reliable scale (and keep the smooth
examples clear of this area), it is performing satisfactorily for its purpose.

The second-to-last test highlights a behavior of the detector that could be
considered a failure mode. Consider a constant function perturbed by white
noise of a much smaller scale. As discussed above, the detector ignores the
constant and only `sees' white noise, yielding a smoothness value of about
zero.  This behavior is undesirable, as the detected smoothness value may
depend on the presence or absence of mere floating point noise.  One root of
this problem is the removal of constant-mode information from the estimation
process, causing the estimator to not have a ``sense of scale'', i.e. keeping
it from noticing that the noise is ``small'' compared to the remainder of the
solution. In the following, we present one way to re-add this ``sense of
scale'' by distributing energy according to a ``perfect modal decay''
\begin{equation}
  \label{eq:viscous-perfect-decay}
  |\hat b_n| \sim \frac 1{\sqrt{\sum_{i=1}^{N_p-1} \frac 1{n^{2N}} }} \frac 1{n^N}
\end{equation}
for $N$ the polynomial degree of the method, where the normalizing
factor ensures that
\[
  \sum_{n=1}^{N_p-1} |\hat b_n|^2 =1.
\]
The idea is to consider the coefficients
\begin{equation}
  \label{eq:viscous-baseline-decay}
  |\tilde q_n|^2\assign
  |\hat q_n|^2 + \|q_N\|^2_{L^2(\D_k)} |\hat b_n|^2
  \qquad \text{for $n\in\{1,\dots,N_p-1\}$}
\end{equation}
as input to skyline pessimization instead of the ``raw'' coefficients $|\hat
q_n|^2$. This amounts to adding \emph{baseline modal decay} scaled by the
element-wise norm that will `drown out' the floating point noise.

For the sake of exposition, baseline decay was not introduced initially.  The
reader may convince himself that its introduction does not unduly modify
experimental results so far by examining the estimated decay exponents given as
``BD+SL'' in the past graphs and comparing to the pure-skyline values given as
``SL''.

\viscousshowestimation{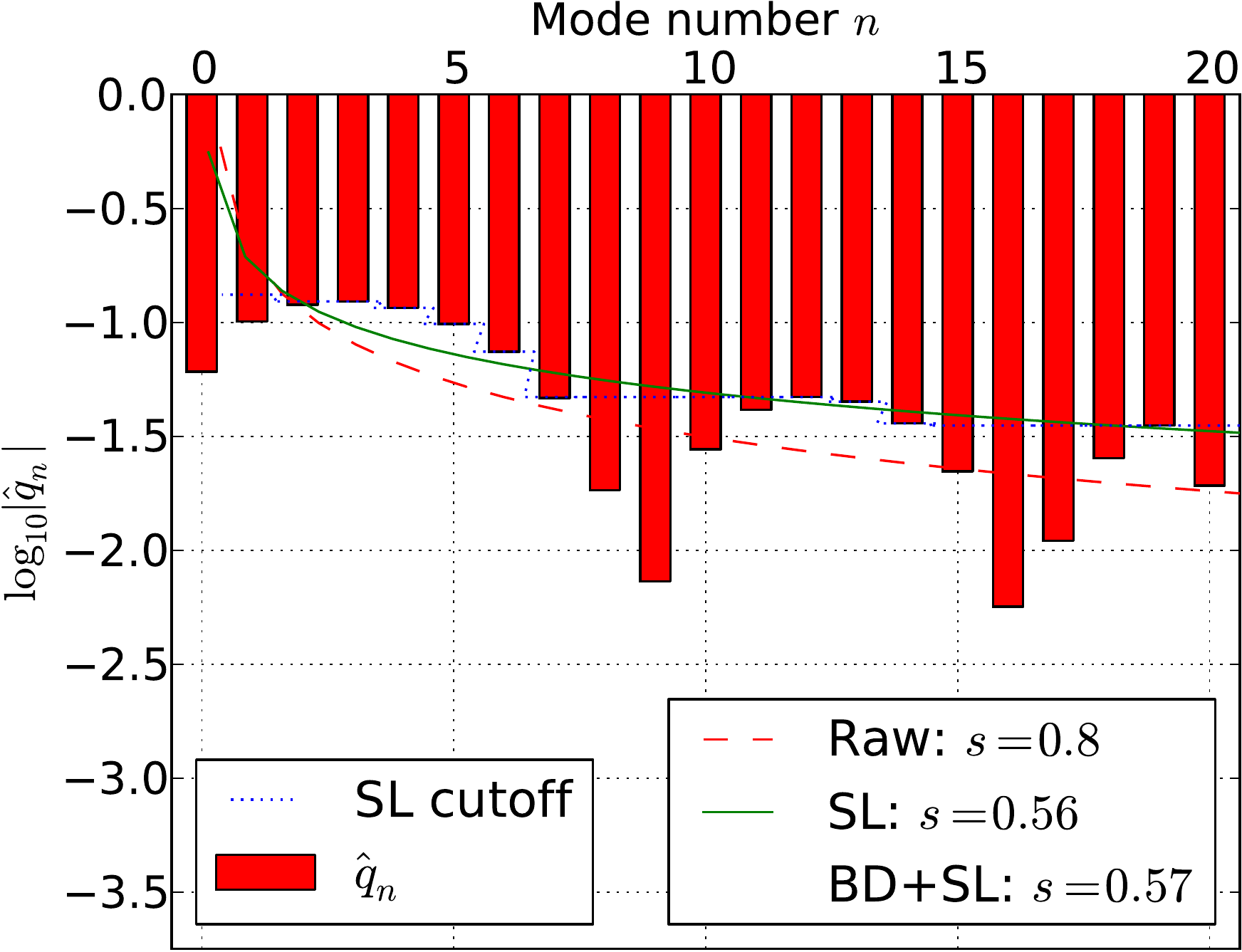}{
  Modal portrait for an approximant of a (discontinuous) jump 
  function, offset from the center of the element.
}

This completes the discussion of the design of the detector. Now might
also be a good time to point out a known shortcoming in its design
that was already anticipated in the motivating discussion. The issue
relates to the discussion of mode scaling with decreasing smoothness
initiated earlier in this section.  Consider Figure
\vref{fig:viscous-estimation-offset-jump}, which shows decay
estimation data for the same Heaviside jump function as Figure
\vref{fig:viscous-estimation-jump}, but shifted to the element's edge.
The data in the figure confirms the earlier conjecture that a function
with a sharply localized non-smoothness might result in modal decay
exponents that differ by up to a factor of two, depending on where the
non-smoothness is located inside the element--the measured smoothness
exponent for the shifted Heaviside function is only 0.57, compared to
1.05 after all corrections above.  Additional confirmation comes from
the fact that the final smoothness estimates for boundary-shifted
versions of the kink and the $C^1$ spline are $s=1.19$ and $s=2.24$
respectively (not shown). This relates in striking ways to the scaling
of the DG CFL condition \eqref{eq:viscous-timestep-restriction}, and
like in its case, an entirely practical remedy for this issue is not
yet known. 

Based on the shown examples, it should be clear that even the
unassisted decay fit is a more robust smoothness estimator than the
single-mode indicator \eqref{eq:viscous-pp-indicator}, if only for the
simple reason that it considers a much broader set of modal data. But
we have shown that even this fairly robust indicator can give poor
results in surprisingly common cases. We feel that this strongly
supports the statement that the decay fit indicator \emph{with}
skyline pessimization and added baseline decay represents a more
practical--if slightly more expensive--way of obtaining smoothness
information on a numerical solution.

\section{From Smoothness to Viscosity}
\label{sec:viscous-smoothness-to-viscosity}
\subsection{Scaling the Viscosity}
\label{sec:viscous-scaling}

This section assumes that the output of the indicator is an estimated
decay exponent $s$, approximating the decay of the solution's modal
coefficients as $|\hat u_n|\sim n^{-s}$. We are seeking to design an
activation function $\nu(s)$ whose value is the viscosity coefficient.

For the interpretation of the decay exponent $s$, recall the targeted
scaling of the smoothness exponent $s$, where (roughly) $s=1$ would
indicate a discontinuous solution, $s=2$ would indicate a $C^0$
solution, $s=3$ a $C^1$ solution, and so forth.  Among the chief
nuisances of polynomial approximations that this work seeks to remedy
is the Gibbs phenomenon, which occurs for discontinuous solutions
($s=1$). We therefore expect to have $\nu(1)=\nu_0$, where
$\nu_0$ is the maximum value of $\nu$ and dictates its
scaling. Merely continuous functions still pose somewhat of a problem
for polynomial approximation, so we arbitrarily fix
$\nu(2)=\nu_0/2$, and finally we fix $\nu(3)=0$, as we
prefer that $C^1$ solutions should not be modified by viscosity.

Given the activation map $\nu_k(s_k)$ of \eqref{eq:viscous-pp-activation-map}
with the fixed values $s_0=2$, the map $\nu(s):=1-\nu_k(s_k)$ with the fixed
values $s_0=2$ and $\kappa=1$ provides such a ramp. (Observe that in
\eqref{eq:viscous-pp-activation-map}, decreasing values indicate more smoothness,
while this work uses the opposite convention.) Because of the close attention
paid to precise scaling of the smoothness $s$, we were able to fix
values for the ramp location and width parameters $\kappa$ and $s_0$.

To find an appropriate value $\nu_0$, the behavior of the diffusion
term needs to be investigated.  To this end, we examine the fundamental
solution of the diffusion equation $u_t=\nu \triangle u$, the \emph{heat
kernel}.  Adopting the probabilistic standard deviation $\sigma$ as a
measure of width, the heat kernel after time $t$ has a width of
$\sigma=\sqrt{2\nu t}$.  Considering some unit $t$ of time, the
conservation law will propagate information to a distance of
$\lambda$, where $\lambda$ is some local characteristic velocity. Observe
that viscosity propagates the bulk of its mass at a non-linear
square-root pace, while the conservation law observes a linear speed.
One therefore needs to pick a reference time scale $t$ as well as a
reference distance at which the two propagation distances
are to coincide.

Choosing $\sigma=h/N$ after $t=(N/2) \Delta t$, and
approximating $\Delta t\approx h/(\lambda N^2)$, one obtains
\begin{equation}
\label{eq:viscosity-scaling}
\nu_0=\frac {\sigma^2} {2t} = \lambda \frac h N.
\end{equation}
This reproduces the value of \citet{barter_shock_2009} and
simultaneously provides some more detailed insight into its meaning.
We would like to note that $\sigma=h/N$ is probably too ambitious a
goal, as this would only smooth discontinuities to a with of about the
distance between two nodal points--likely too little as Figure
\vref{fig:viscous-estimation-jump} shows.  A choice of $\sigma=3h/N$
has proven to be more realistic.

For a system of conservation laws, there remains the question of which
characteristic velocity should be chosen for $\lambda$. This choice
has important implications as, e.g. in the Euler system, contact
discontinuities propagate with stream velocity, whereas shocks
propagate at sonic speeds.  In a one-dimensional setting,
\citet{rieper_dissipation_2010} convincingly argues that the best
course of action is to perform smoothing in characteristic variables,
so that each wave receives the amount of smoothing specified by the
scheme, e.g. as given in \eqref{eq:viscosity-scaling}. Observe that
doing may work well in one-dimension and for low-order multi-D finite
volume schemes, but it is less clear how it might be applied in a
genuinely multidimensional situation.  A simple and functional
strategy is to choose $\lambda$ to be the maximum characteristic
velocity $\lambda_{\text{max}}$. The simplicity of this strategy comes
at a price, however: returning to the example of the Euler equations,
contact discontinuities have their $\nu_0$ set higher than would be
necessary from this analysis, and our numerical experiments will
reflect this.

Thus the $\lambda_{\max}$-based scaling is not perfect. It works,
in the sense that all test examples run successfully using it, but
some can benefit from an additional `fudge factor'. For example, while
problems involving Burgers' equation (not shown) work well
with an unmodified scaling in a 'picture norm' sense (little
oscillation, least smoothing), most subsonic Euler problems benefit
from the application of an additional factor of $1/2$.  This is not
entirely unexpected, given the above discussion.

\subsection{Smoothing the Viscosity}
\label{sec:viscous-smoothing}

The artificial viscosity $\nu(x)$ obtained so far is a per-element
quantity, with no guarantees on how it might vary across the domain.
In particular, since the viscosity is constant on each element, it
will invariably be discontinuous.

Now observe how the viscosity is employed in the equations of Section
\ref{sec:viscous-equations}. In particular, observe that in order to
maintain conservativity, the viscosity occurs \emph{inside} a
derivative. Great care is required in the correct numerical solution
of a diffusion equation with discontinuous viscosities using
discontinuous Galerkin methods.  \citet{proft_discontinuous_2009,
loercher_explicit_2008, ern_discontinuous_2009} describe various
precautions that need to be taken to avoid non-conservativity and
non-consistency.

\citet{feistauer_robust_2007} also notice the issues caused by
localized, discontinuous viscosities and  propose an adapted flux term
to ``strengthen the influence of neighbouring elements and [improve]
the behaviour of the method''. \citet{barter_shock_2009}, through
numerical experiment, also arrive at the conclusion that a discontinuous
viscosity causes issues and show a marked decrease in $H^1$ error for
smooth viscosities.  Since one is at considerable liberty to choose
the viscosity $\nu(x)$, we agree that it is best to choose a $\nu$ that
does not include discontinuities, to avoid this entire complex of
issues.

Therefore, given that the detection infrastructure built up so far
works in an element-by-element fashion, one needs to introduce a
post-processing step that somehow smoothes out generated $\nu$. In
doing so, one again has a wide array of choices.
\citet{barter_shock_2009} propose a diffusion equation (effectively
``diffusing the diffusivity'') with time-relaxation to obtain a
viscosity that is smooth in both time and space. Unfortunately, this
choice is unsuitable given the design choices for explicit time
stepping laid out in Section
\ref{sec:viscous-design-considerations}--to achieve sufficient
smoothing of the viscosity, one needs to choose a large diffusivity
for it, which results in a very stiff system of ODEs.

One important question in the design of a successful smoothing method
is, precisely how smooth must the result of the smoothing be? In
computational experiments relating to artificial
viscosity, we have found that there does not seem to be an advantage to
having the viscosity $\nu\in C^k$ for $k>0$. 

Based on these considerations, the method employed in the experiments in the
next section proceeds as follows:
\begin{enumerate}
  \item At each vertex, collect the maximum viscosity occurring in
  each of the adjacent elements.
  \item Propagate the resulting maxima back to each element adjoining
  the vertex.
  \item Use a linear ($P^1$) interpolant to extend the values at the
  vertices into a viscosity on the entire element.
\end{enumerate}
In our experience, this method is cheap, reasonably straightforward to
implement even on GPUs, and it satisfies the design requirements set forth
above.

\section{Experience with and Evaluation of the Scheme}
\label{sec:viscous-experiment}

\newcommand{\viscsubfigpic}[2]{
  \subfigure[#2]{
    \includegraphics[width=0.4\textwidth]{#1.pdf}
    \label{fig:viscous-#1}
  }
}

\subsection{Advection: Basic Functionality, Interaction with Time
Discretization}

\begin{figure}
  \begin{center}
  \viscsubfigpic{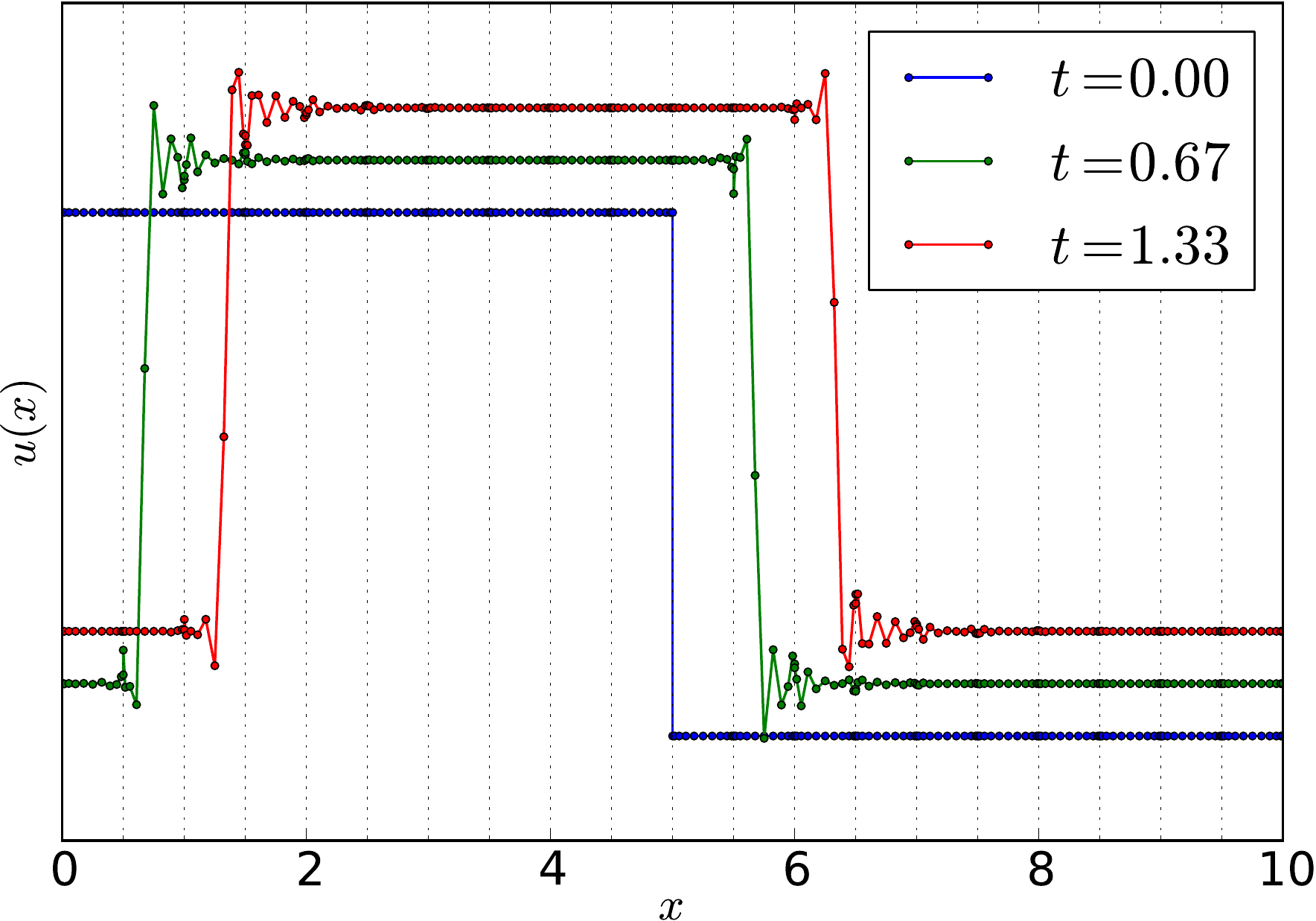}{
    Solution of the advection equation without artificial viscosity.
  }
  \viscsubfigpic{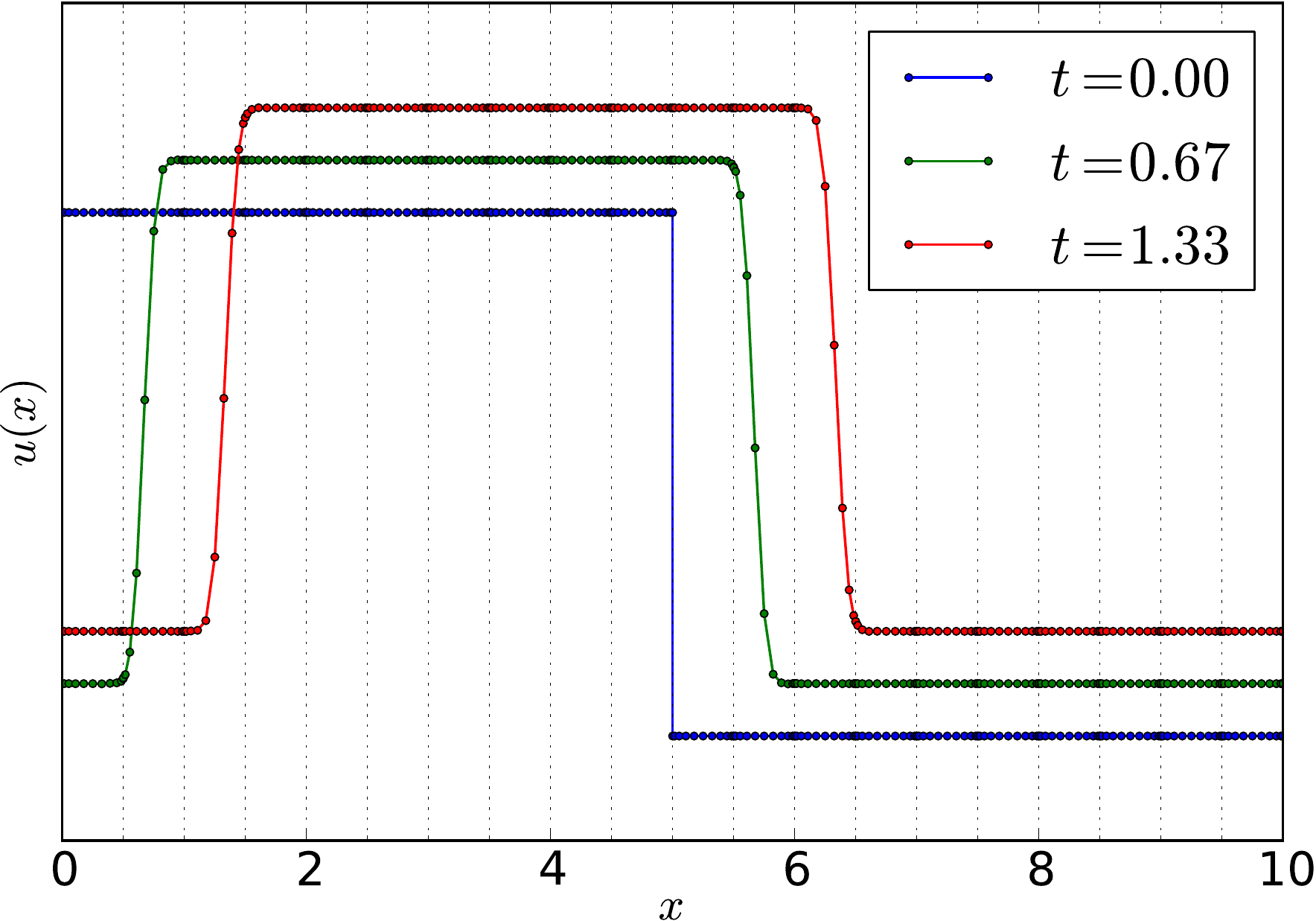}{
    Solution of the advection equation with artificial viscosity,
    after short amounts of time.
  }
  \end{center}
  \caption{Spatial shock capturing behavior of the artificial
  viscosity scheme on an advection equation.}
\end{figure}

The first set of results we would like to discuss relates to the
advection equation (Section \ref{sec:viscous-eqn-advection}).
The examples in this section examine the advection of the function
$u(x):=\mathbf 1_{[0,5)}$ over an interval $(0,10)$.

\citet{kuzmin_flux_2005} suggest that the advection equation is
particularly suited to testing shock capturing schemes for two
reasons: First, because it is the simplest PDE that can sustain a
discontinuous solution, so that the behavior of the method can be
observed in a well-understood setting, isolated from other
characteristics and nonlinear effects.  Second, because
discontinuities in it are not self-steepening, in analogy to contact
discontinuities in the Euler equations, it makes a challenging example
to be treated with artificial viscosity: Once a discontinuity is
unduly smeared by viscosity, nothing will return it to its former,
sharp shape.

Figure \vref{fig:viscous-advection-without-visc} displays the behavior
of the unmodified discontinuous Galerkin method as described in
Section \ref{sec:viscous-eqn-advection}. As expected, a strong
Gibbs-type overshoot is observed, although it is worth noting that the
used upwind fluxes already provide enough dissipation of
high-frequency modes to prevent the solution from becoming useless.
This example, and all examples that follow in this subsection, were
run at polynomial degree $N=10$ on a discretization using $K=20$
elements.

Next, Figure \vref{fig:viscous-advection-with-visc} displays the
result of the same calculation once the artificial viscosity machinery
as described above is enabled. Discontinuities are resolved within
eight points, i.e. within less than one element (containing $N_p=11$
points) and have no visible overshoots.  (Note that as an expected
consequence of the clustering of the nodes towards element edges,
points appear spaced closer together where the discontinuity touches
an element boundary.) Element boundaries are shown as dashed lines for
orientation.  Figure \ref{fig:viscous-advection-with-visc} displays
the solution after only a brief amount of simulation time has passed.
It turns out that the solution--at least visually--settles into its
final form and does not change much even after a large number of
round-trips. The steepness of the solution is retained as in Figure
\ref{fig:viscous-advection-with-visc}, and the number of points that
are required to resolve the discontinuity remains stable.

\begin{figure}
  \begin{center}
  \viscsubfigpic{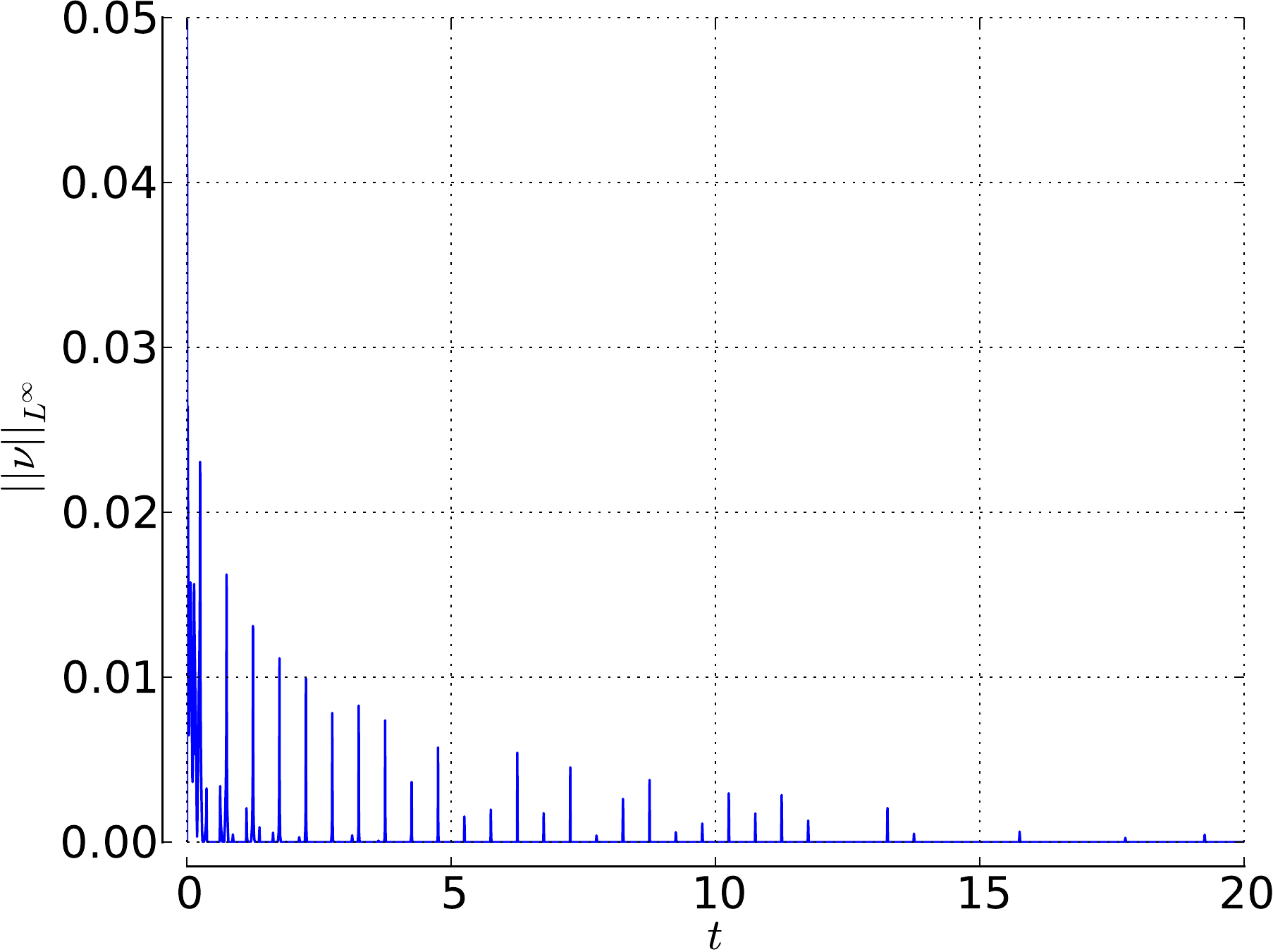}{
    Artificial viscosity activations vs. time, in a discontinuous
    advection calculation.
  }
  \viscsubfigpic{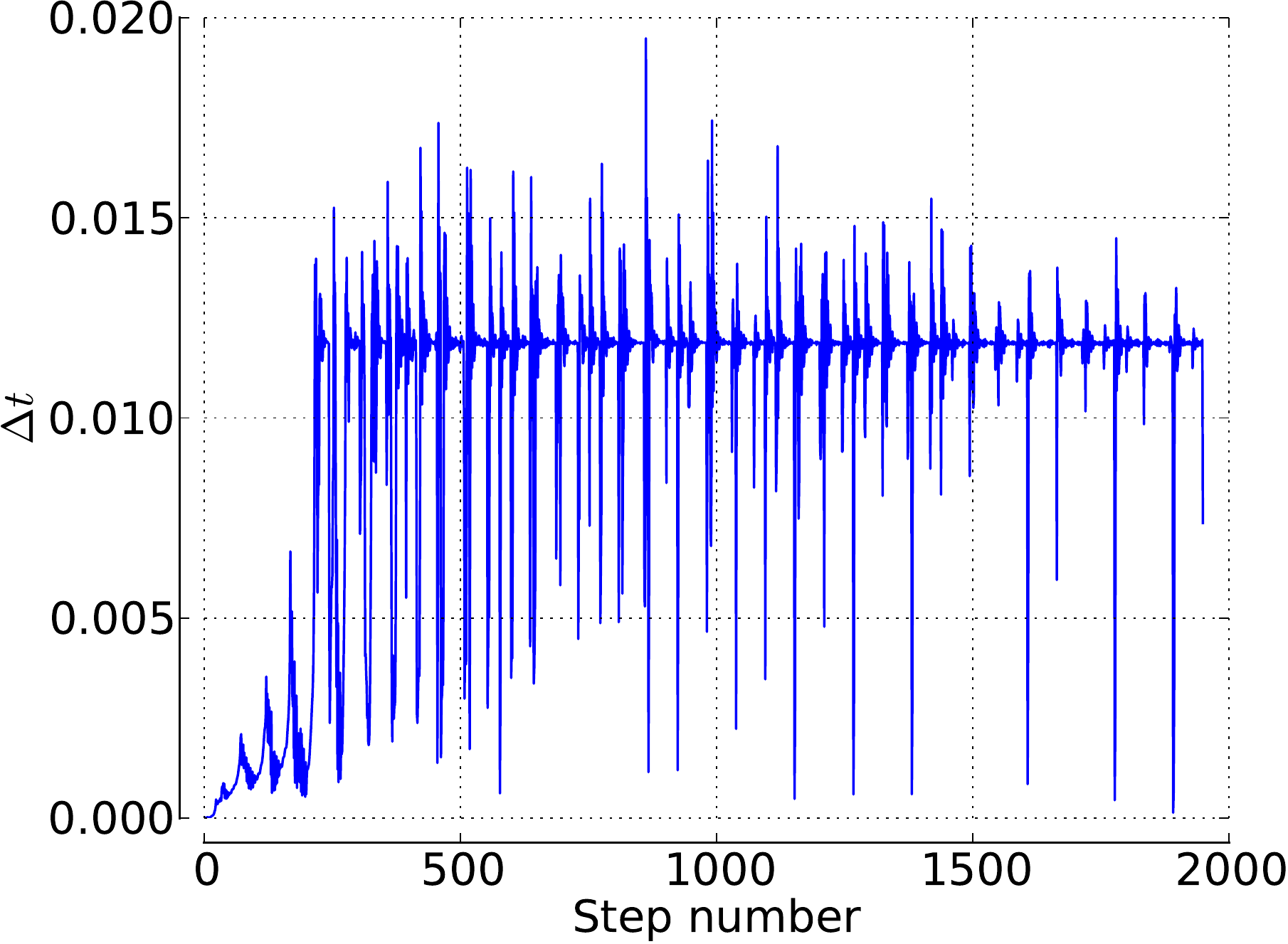}{
    Adaptively found time step vs. step number, in a discontinuous
    advection calculation.
  }
  \end{center}
  \caption{Interaction of the shock-capturing artificial viscosity
  with the time discretization.}
\end{figure}

Figure \vref{fig:viscous-advection-element-transitions} sheds a new
light on this ``settling'' observation and the observed increased
sensitivity of the detector near element boundaries that was discussed
above.  It shows the maximum viscosity $\|\nu\|_{L^\infty}$ found
anywhere on the domain, graphed versus simulation time. If the
observation of ``brief-settling-then-steady-state'' were entirely
true, then one would observe no sensor activations whatsoever after
``settling'' has occurred. This is not what is observed here. Instead,
one sees a slowly decaying train of viscosity activation spikes. It
turns out that each of these spikes coincides with a discontinuity
crossing an element boundary. This again confirms the observation that
the detection scheme is inhomogeneous in space, i.e. it judges
solution smoothness differently depending on whether a discontinuity
is located in the interior of an element or at its boundary. Since the
sensor is only exposed to the non-smoothness for very short periods at
a time, according to Figure
\ref{fig:viscous-advection-element-transitions} it takes considerable
time ($t\gtrapprox 12$ in the example) and a number of viscosity
``spikes'' until a profile is achieved that does not trip even the
overly sensitive version of the detector. It is to be expected that
the final profile is twice smoother than would be required if the
oversensitivity did not exist.

As a last observation on the behavior of the method on this
exceedingly simple problem, we would like to examine its interaction
with the adaptive time stepper. The examples were computed using the
well-known embedded Runge-Kutta method of third order by
\citet{bogacki_pair_1989} (``\texttt{ode23}'' in Matlab).
\vref{fig:viscous-advection-delta-t} shows the adaptively-chosen time
step $\Delta t$ as a function of the step number. The stable advective
time step is clearly visible, as is the initial ``settling'' period
discussed above, along with a variety of time step reductions
occurring along the way. Some of these coincide with element
transitions of discontinuities, but the situation is more ambiguous
(and noisier) than in the case of viscosity activations. The figure
does make one thing amply clear, however: an
artificial-viscosity-based shock capturing scheme using explicit time
stepping must use time step adaptivity, or it will not be competitive.

\subsection{Waves: Shock Spreading and Spurious Coupling}
\label{sec:viscous-wave-observations}
The next, more complicated problem for which we examine the behavior of
the proposed artificial viscosity is the wave equation, described in
Section \ref{sec:viscous-eqn-wave}.

\begin{figure}
  \centering%
  \viscsubfigpic{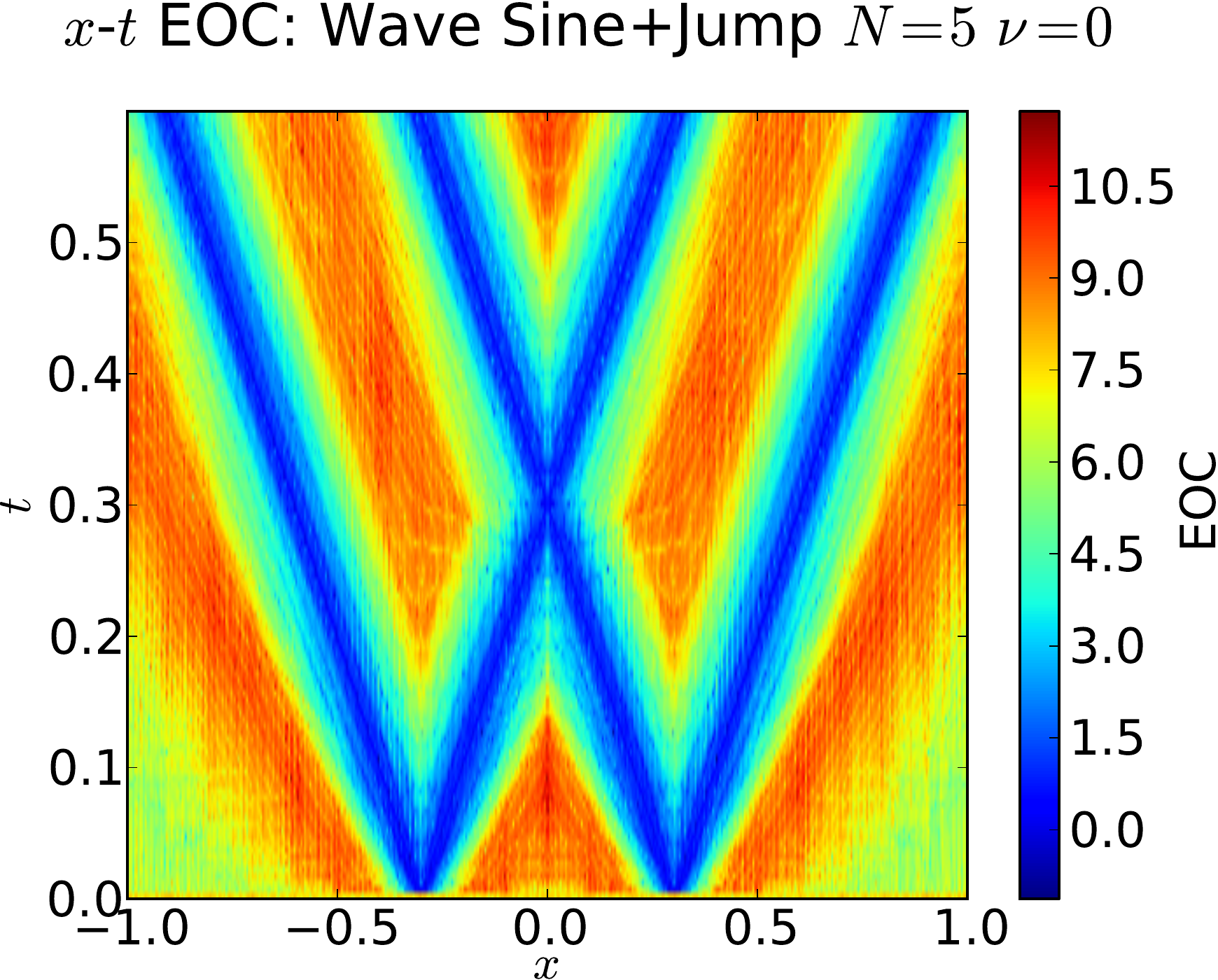}{
    EOC for the wave equation with a discontinuous initial condition
    \emph{without} artificial viscosity.
  }%
  \viscsubfigpic{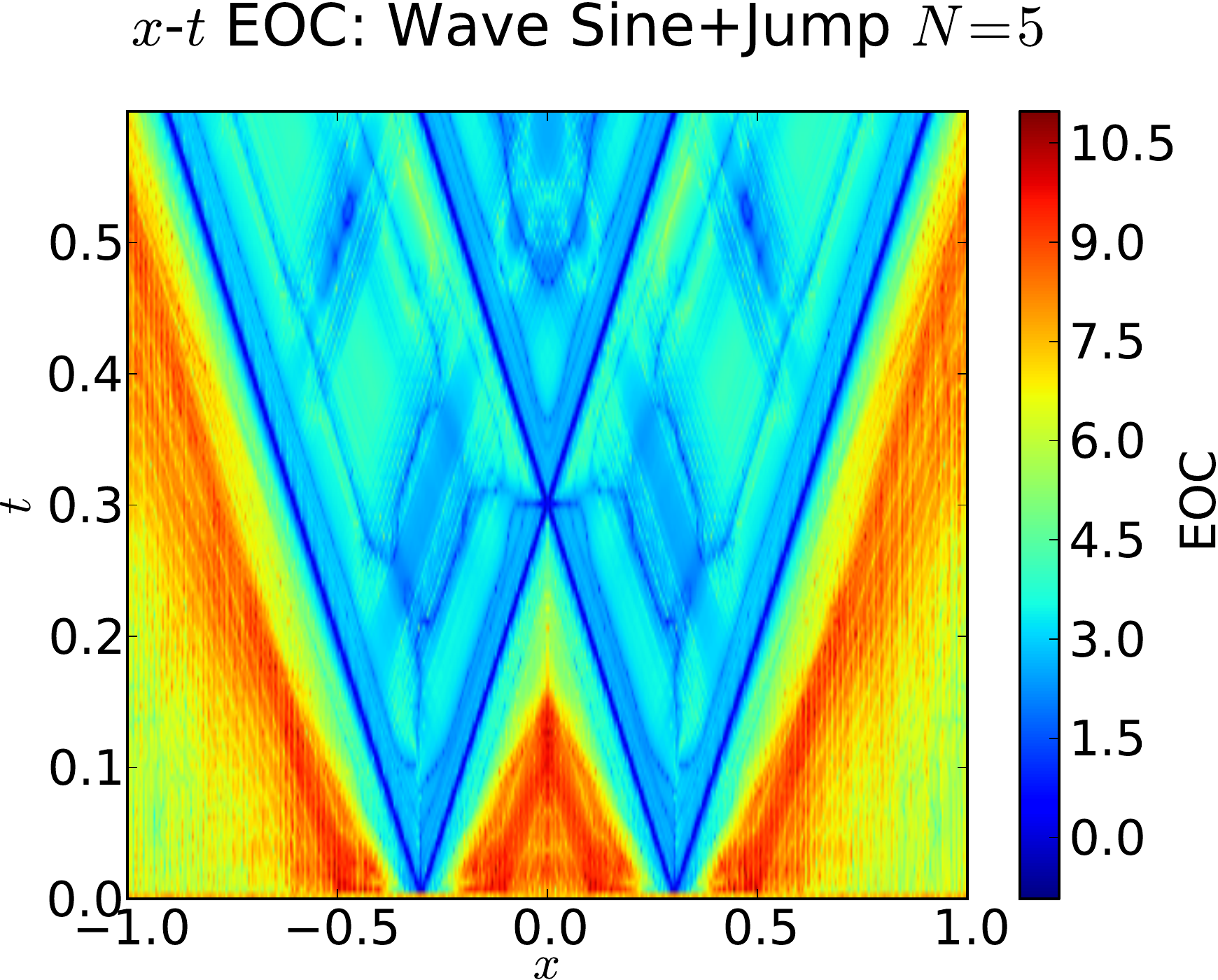}{
    EOC for the wave equation with a discontinuous initial condition
    \emph{with} artificial viscosity.
  }%
  \caption{
    Empirical order of convergence for the wave equation with
    discontinuous initial conditions.
  }
  \label{fig:viscous-spatial-eoc-wave}
\end{figure}

We would like to set the stage for our experimental results by
considering the context of recent work by \citet{cockburn_error_2008},
who show (under a number of additional assumptions) that for a DG
computation of a linear advection equation at second order using a
second-order total-variation-diminishing (TVD) time discretization,
pollution of the numerical solution by the shock by time $T$ stays
localized to an area of size $O(\sqrt{hT})$ ahead of and an area of
size $O(\sqrt[3]{Th^2})$ behind the discontinuity.  Although they only
show this for a scalar advection equation, the wave equation
\eqref{eq:viscous-wave-eqn} and its discretization may be transformed
into two decoupled advection equations, and hence the result applies
in this case as well.

We will study the pollution of the solution by examining its pointwise
empirical order of convergence to the known analytic solution in space
and time, starting from the initial condition
\[
  u(x,0) = 2+\cos(5\pi x) + 4\cdot \mathbf{1}_{[-0.3,0.3]}(x),
  \qquad
  v(x,0) = 0,
\]
subject to Neumann boundary conditions, on a domain $\Omega=(-1,1)$ up
to a final time $T=0.6$, with a wave speed $c=1$.

Figure \vref{fig:viscous-spatial-eoc-wave} shows the resulting
convergence plots, obtained with and without artificial viscosity. As
expected through the work of \citet{cockburn_error_2008}, the inviscid
DG scheme of Figure \ref{fig:viscous-wave-no-visc-xt-eoc} achieves
full convergence away from the discontinuities, but also shows a
slowly-growing zone of non-convergence near the discontinuities, again
matching predictions.

Unfortunately, results are not as favorable once artificial viscosity
starts to act on the scheme. Outside the region that interacts with the
discontinuities, convergence is roughly as before. However inside the
interacting regions, convergence does improve again away from the
discontinuity, but it does not recover the full order of the scheme.
This reduction in order is in line with results obtained for
finite-difference solutions downstream of a slightly viscous shock by
\citet{efraimsson_remark_1999} (see also
\citep{kreiss_elimination_2001}). The observation further underscores
the importance of the wave equation as a test example for shock
capturing schemes.  Once the PDE is rewritten in as a system of
first-order conservation laws, the single
added viscosity of \eqref{eq:viscous-wave-eqn} induces a
cross-coupling that appears to destroy accuracy.

Note that such behavior \emph{cannot} be observed in the advection
equation, or, generally, any purely scalar conservation law, since
these equations have only one characteristic wave, and hence the
pollution caused by the artificial viscosity cannot spread, but
propagates along with the solution.  This might lead one to suggest
an obvious ``fix'' for the issue: The first-order system (i.e.
the left-hand side of \ref{eq:viscous-wave-eqn}) can
easily be transformed into characteristic variables, where it takes
the form of two advection equations that only couple at the boundary,
such that the issue disappears \citep{rieper_dissipation_2010}.  As we
have already discussed, proposing this is as a general remedy is
however a bit disingenuous, as it cannot work properly in multiple
dimensions. Another idea that one might have to try and avoid the
reduction in accuracy is to use separate viscosities for each of the
variables. According to our experiments, this does not help, as the
cross-coupling of the system persists.

Next, it seems unlikely that this problem is specific to the
artificial viscosity constructed in this article, or to discontinuous
Galerkin methods, for that matter. It should be
investigated whether \emph{all} artificial viscosity schemes proposed
so far in the literature suffer from this shortcoming.

\subsection{Euler's Equations}
In this section, we will carefully examine the behavior of the
artificial viscosity method introduced above on Euler's equations of
gas dynamics, starting with the classical exact solution of the
Riemann problem given by \citet{sod_survey_1978} as the first example.

\begin{figure}
  \centering%
  \viscsubfigpic{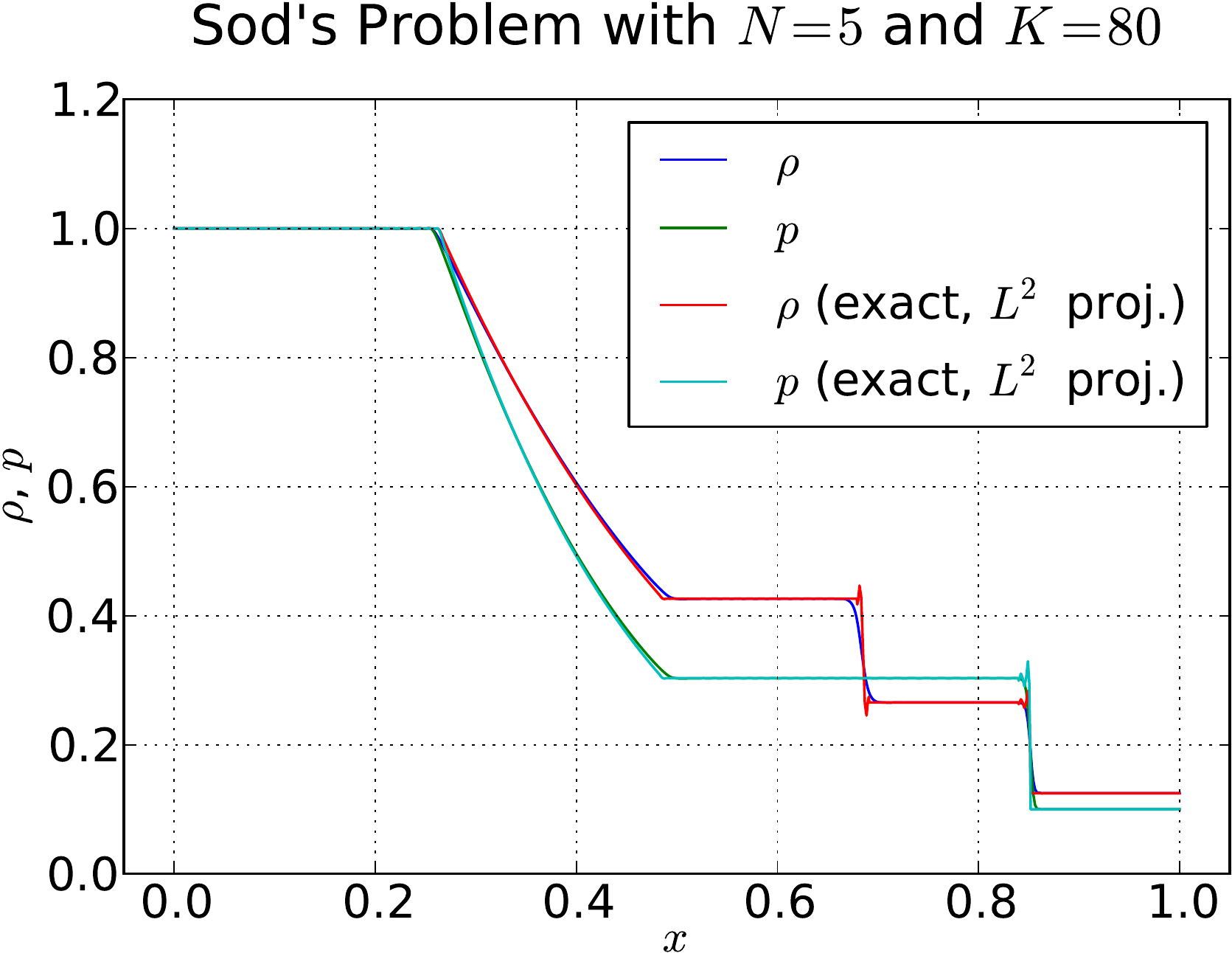}{
    $L^2$-projected exact and approximate numerical solutions of Sod's
    problem for polynomial degree $N=5$ in $K=80$ elements.
  }%
  \viscsubfigpic{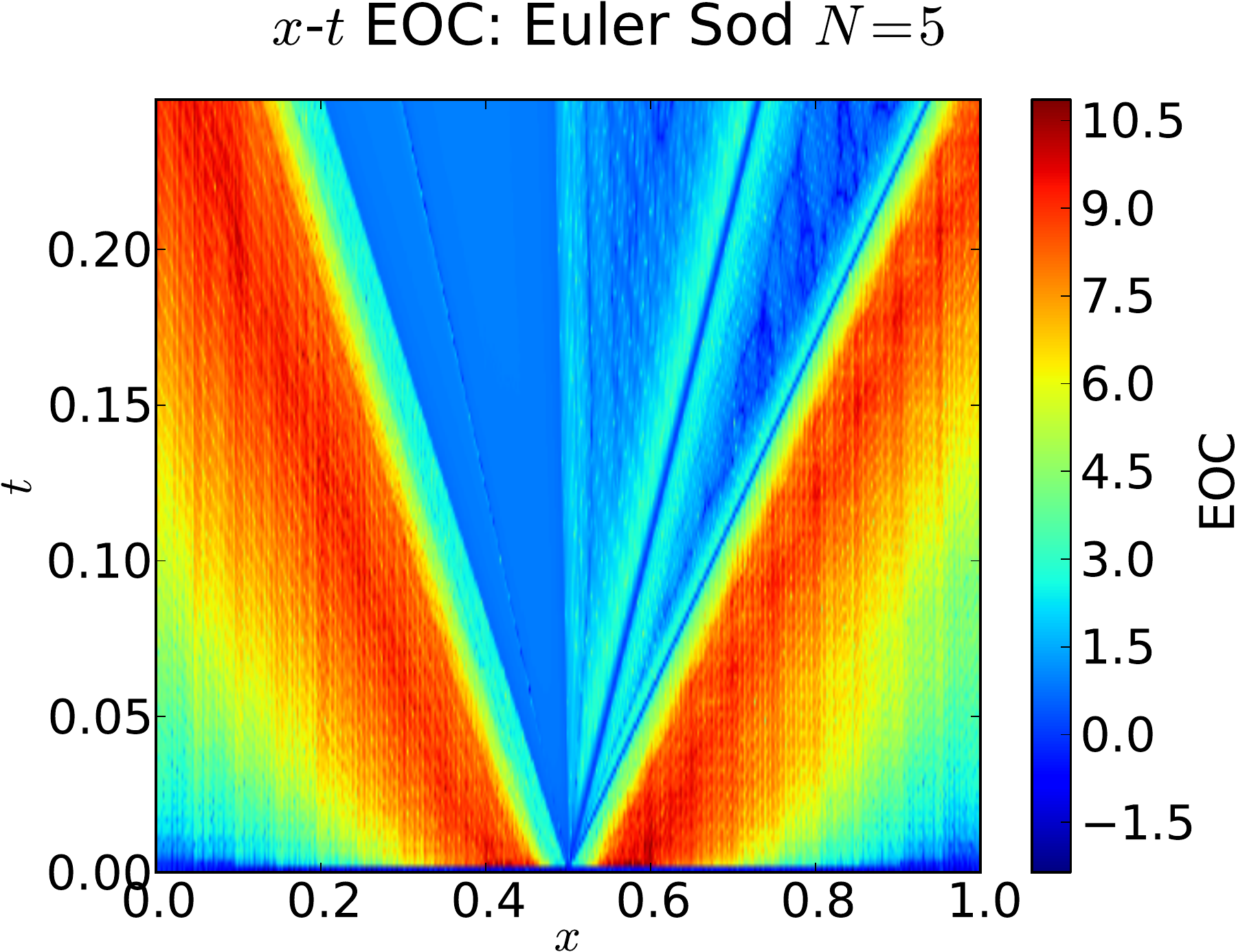}{
    Space-time diagram of the empirical order of convergence for Sod's problem,
    computed with artificial viscosity.
  }%
  \caption{Sod's problem with artificial viscosity: solution and
  $x$-$t$ convergence. }
  \label{fig:viscous-sod}
\end{figure}

Figure \vref{fig:viscous-sod-n5-k80} shows computational results,
again at polynomial degree $N=5$ on $K=80$ elements, in direct
comparison with the ($L^2$ projection of) the exact solution, for the
density $\rho$ and the pressure $p$, at the final time $T=0.25$ of the
computation.

While the figure above gives an impression of the desired solution and
a first impression of the performance of the method, it is perhaps
more enlightening to examine an analog to the the convergence in space
and time of Figure
\vref{fig:viscous-spatial-eoc-wave} in the gas dynamics setting.  Figure
\vref{fig:viscous-euler-sod-xt-eoc} provides this. As above, the
computation was carried out at polynomial degree $N=5$, at a variety
of mesh resolutions ranging from $K=20$ to $320$ elements across the
domain.  Like in the linear case, convergence away from the shock
region is good, while in the central, shock-interacting `fan', it hardly
exceeds order 1. In particular, it is worth noting that convergence
along the profile of the smooth rarefaction wave is also no better
than order 1. Given the results obtained for the wave equation, this
is not very surprising, and it confirms that the issues observed on
linear problems persist in the nonlinear case.

\begin{figure}
  \centering
  \viscsubfigpic{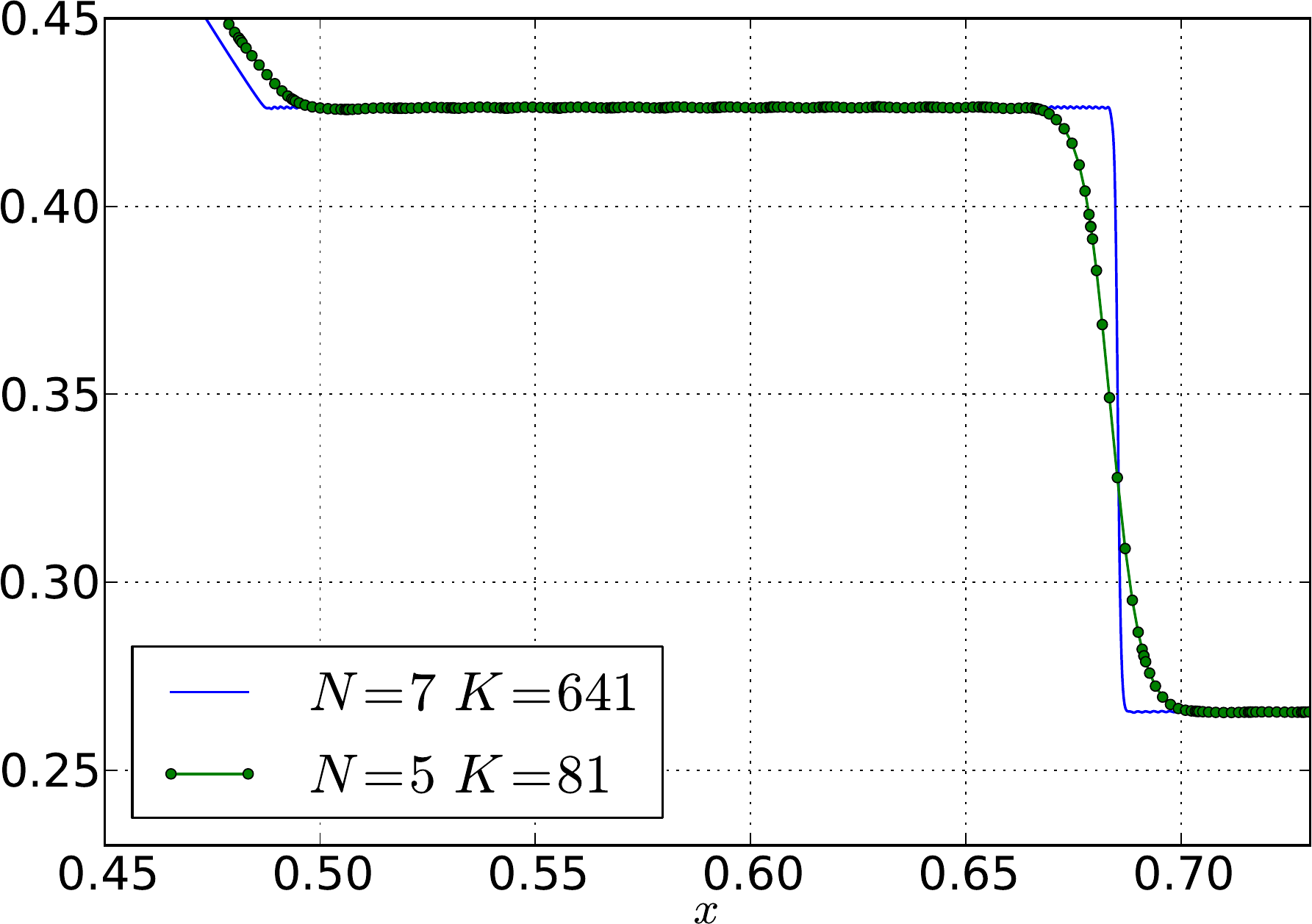}{
    Close-up view of the contact discontinuity in Figure
    \vref{fig:viscous-sod-n5-k80} at low and high numerical
    resolutions. Interpolation nodes for the low-resolution case are
    shown as dots.
  }
  \viscsubfigpic{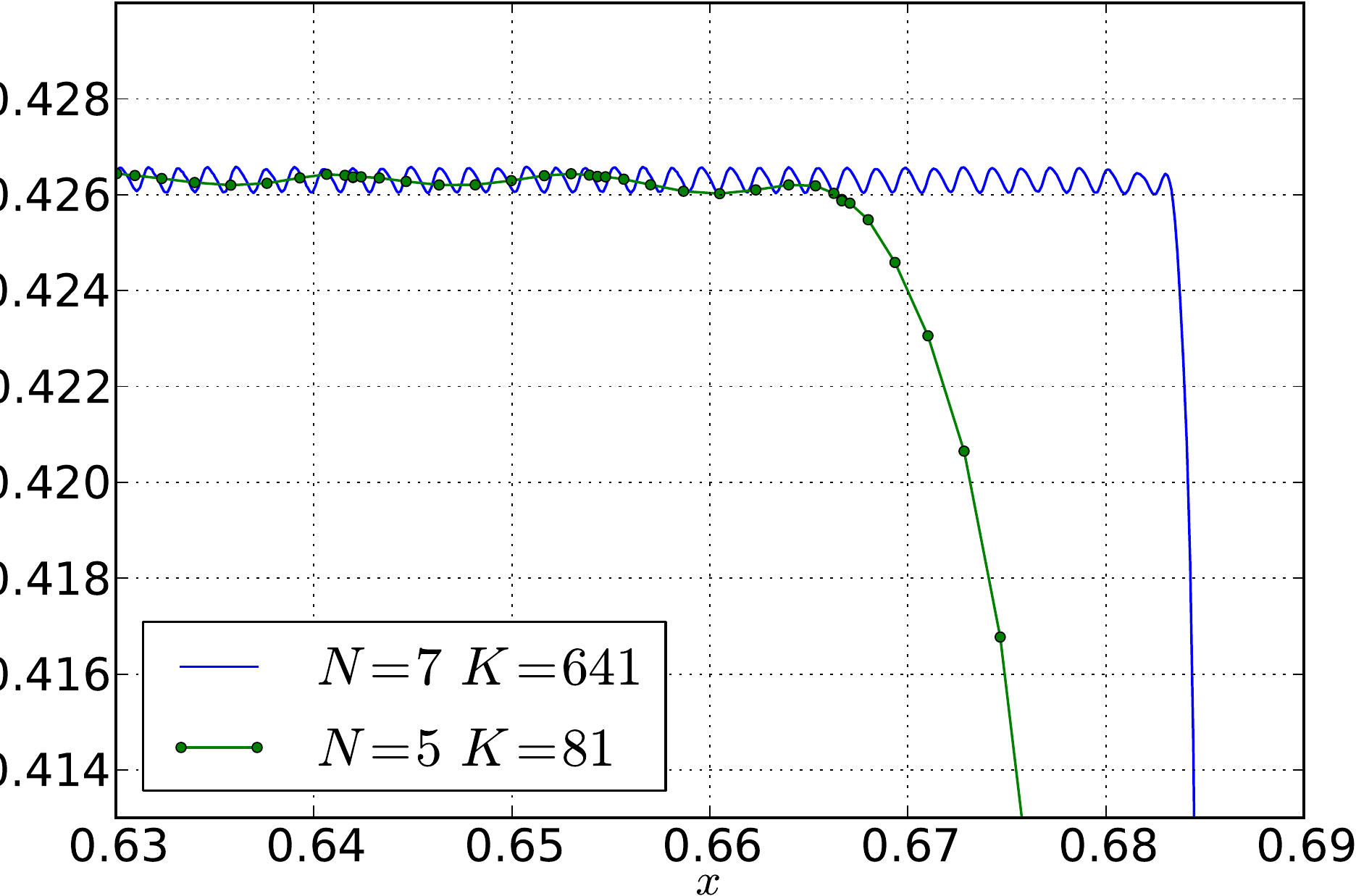}{
    Extreme close-up view of the tip of the contact discontinuity in
    Figure \vref{fig:viscous-shock-wrinkles-cd}, at low and high
    numerical resolutions.
  }
  \caption{Element-scale oscillation exhibited by the artificial
    viscosity scheme.}
  \label{fig:viscous-element-scale-oscillation}
\end{figure}

A closer look at the numerical solutions in the poorly-converged
region of \ref{fig:viscous-euler-sod-xt-eoc} offers a revealing
insight, shown in Figure \vref{fig:viscous-element-scale-oscillation}
for a high-resolution case ($N=7$, $K=641$) and a low-resolution case
($N=5$, $K=81$).  On the constant parts of the solution to the Riemann
problem, we observe small ``wrinkles''. Figure
\ref{fig:viscous-shock-wrinkles-cd} provides a sense of scale, while
the extreme close-up of Figure \ref{fig:viscous-shock-wrinkles-big}
shows the phenomenon in detail. In both the high- and the
low-resolution case, the oscillation's wave length roughly agrees with
the size of an element. Further, it is remarkable that the magnitude
of the oscillation appears to grow, rather than shrink, with increased
resolution, which seems to indicate that convergence below the margin
provided for by the oscillation might not occur.  (Convergence will be
examined in some detail below.) The phenomenon is observed on all
constant areas that are inside the fan of characteristics emanating
from the shock at time $t=0$. So far, we do not understand the cause of
this phenomenon, nor is it known whether there is a connection between
these wrinkles and the reduced convergence observed in Section
\ref{sec:viscous-wave-observations}.  One might speculate that, again,
the detector's spatial inhomogeneity is to blame. While we are as yet
unsure of the source of the phenomenon, we would like to  note that post-shock
oscillations of this nature have been observed and studied even in
schemes that do not use element-based decompositions \citep{arora_postshock_1997}.

\begin{table}
  \centering
  \begin{tabular}{l|cccc|c}
   & $N=4$ & $N=5$ & $N=7$ & $N=9$ & EOC \\
  \hline
  $h/1$ & $9.982\cdot 10^{-3}$ & $7.934\cdot 10^{-3}$ & $6.522\cdot 10^{-3}$ & $5.567\cdot 10^{-3}$ & 0.70 \\
  $h/2$ & $5.442\cdot 10^{-3}$ & $4.231\cdot 10^{-3}$ & $3.395\cdot 10^{-3}$ & $2.921\cdot 10^{-3}$ & 0.75 \\
  $h/4$ & $2.945\cdot 10^{-3}$ & $2.219\cdot 10^{-3}$ & $1.778\cdot 10^{-3}$ & $1.568\cdot 10^{-3}$ & 0.76 \\
  $h/8$ & $1.548\cdot 10^{-3}$ & $1.166\cdot 10^{-3}$ & $9.488\cdot 10^{-4}$ & $8.329\cdot 10^{-4}$ & 0.74 \\
  $h/16$ & $8.087\cdot 10^{-4}$ & $6.006\cdot 10^{-4}$ & $5.121\cdot 10^{-4}$ & $4.598\cdot 10^{-4}$ & 0.66 \\
  $h/32$ & $4.207\cdot 10^{-4}$ & $3.111\cdot 10^{-4}$ & $2.806\cdot 10^{-4}$ & --- & 0.69 \\
  \hline
  EOC & 0.93 & 0.95 & 0.92 & 0.92 &  \\
  \end{tabular}
  \caption[
    $L^1$ error and convergence data for the Sod problem of the Euler
    equations of gas dynamics.
  ]{
    $L^1$ error and convergence data for the Sod problem of the Euler
    equations of gas dynamics.
    ``EOC'' stands for the empirical order of convergence, obtained as a least-squares
    fit to the data.
  }
  \label{tab:viscous-euler-sod-l1-convergence}
\end{table}

Beyond the spot testing conducted so far, we have also carried out a
more comprehensive convergence study on the Euler equations applied to
the Sod problem. The raw $L^1$ error data as well as empirical
convergence order results obtained from least-squares fits are shown
in Table \vref{tab:viscous-euler-sod-l1-convergence}. The data was
gathered at a variety of polynomial degrees $N$ and with $K=20$
elements at the coarsest level, with uniform refinements thereafter.
The data seems to support about a full order of convergence in
$h=1/K$. No improvement in convergence occurs as the order is
increased. Further, the data supports less than a full order of
convergence in $N$, indicating that an addition of elemental
resolution at present is a more effective way of getting a more
accurate solution than increasing the size of the local approximation
spaces, especially considering that the computational complexity grows
superlinearly in $N$. At the resolutions examined, the influence of
the oscillations (``wrinkles'') observed above does not appear to have
contributed a significant part of the error--given their observed
behavior in response to resolution changes, they would likely have
represented a ``bottom'' to convergence at some fixed error magnitude.
That issue aside, the observed convergence data appears to be as good
as one might reasonably expect. While convergence of higher order
would of course be desirable, the method as it presently stands is not
designed to be able to achieve this. Through some experiments on
polynomials, we have reason to believe that convergence of order one
in $N$ is achievable and thereby a goal for future research.

In addition to the problem of \citet{sod_survey_1978}, which has
furnished the basis for all tests so far, we have also conducted tests
using other available solutions for the Euler equations. One such
solution that is rather similar to the Sod problem is that of
\citet{lax_weak_1954} in that it also originates from a Riemann
problem. Figure \vref{fig:viscous-lax-n5-k80} demonstrates that the
scheme can successfully compute a correct solution to the problem.
Lax's problem prominently features a contact discontinuity, which is
prone to smearing, as was discussed above.  The contact discontinuity
in the figure appears somewhat more smeared than the Sod contact
discontinuity at a similar scale.

\begin{figure}
  \viscsubfigpic{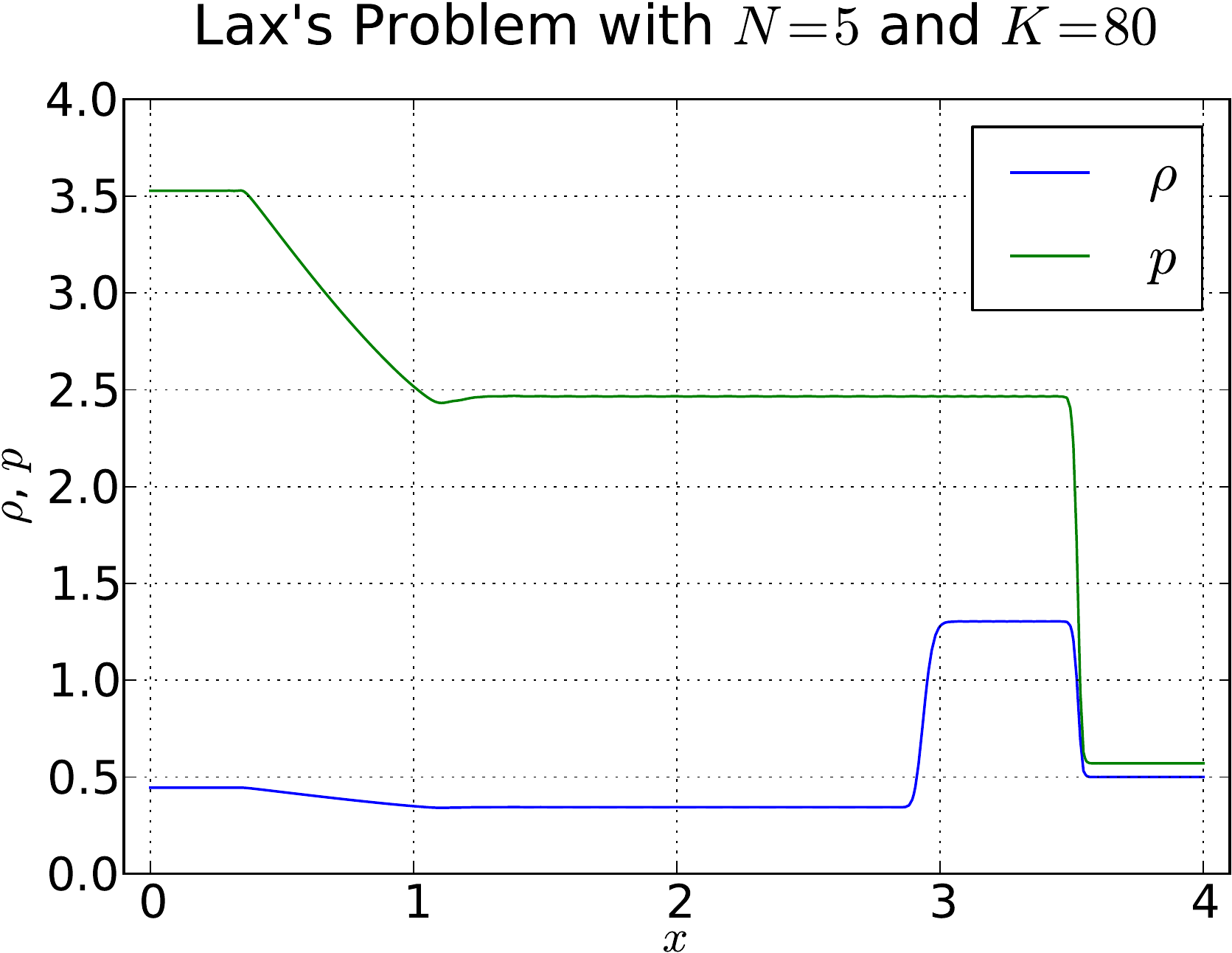}{
    Approximate numerical solution for density and pressure of Lax's
    problem for polynomial degree $N=5$ in $K=80$ elements.
  }
  \hfill
  \viscsubfigpic{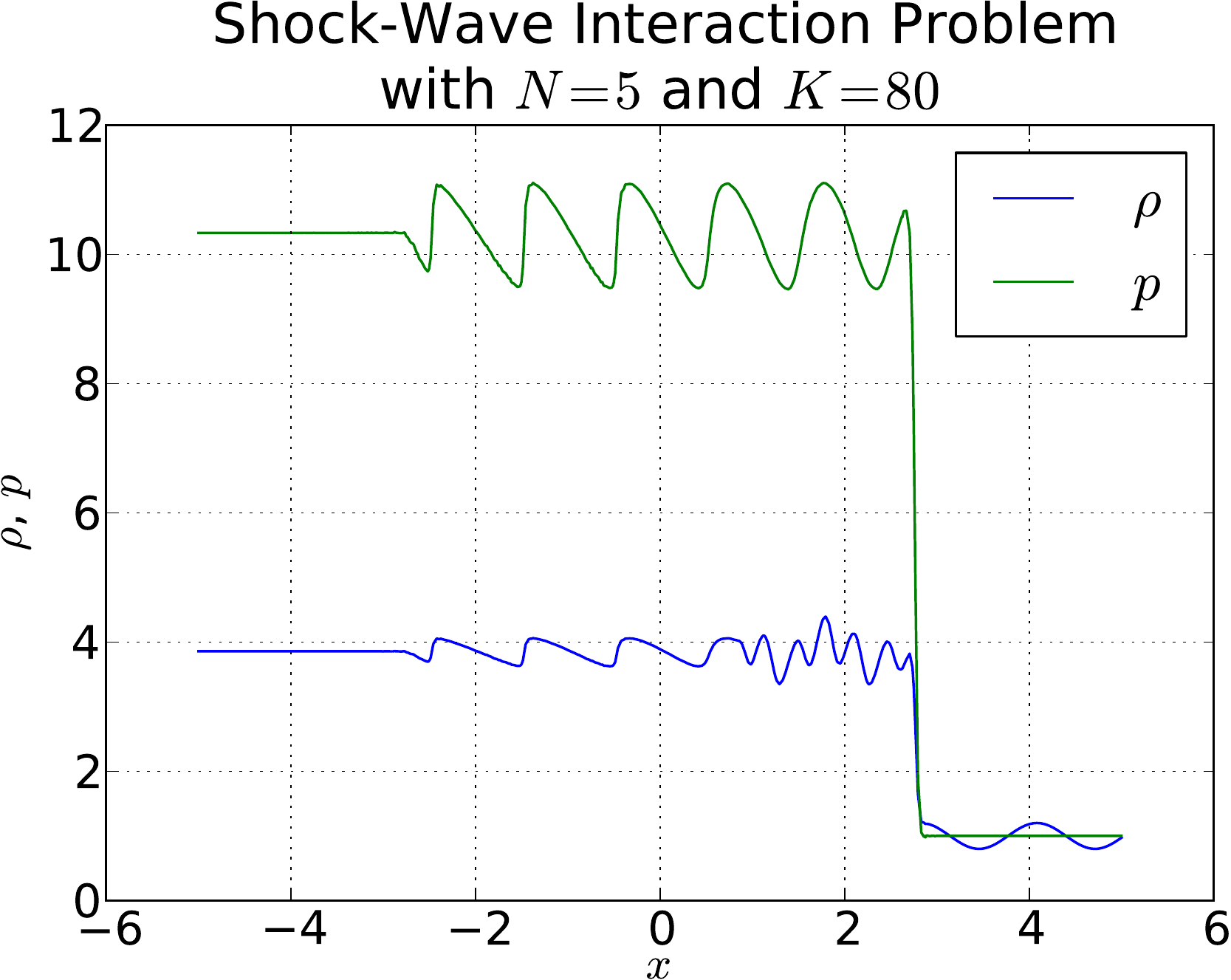}{
    Approximate numerical solution for density and pressure of a
    shock-wave interaction problem for polynomial degree $N=5$ in
    $K=80$ elements.
  }
  \caption{Solutions of classical test problems for the Euler
    equations using the artificial viscosity scheme.}
  \label{fig:viscous-euler-1d-tests}
\end{figure}

A further basic benchmark test for the method applied to the
one-dimensional Euler equations was proposed by \citet[Example
8]{shu_efficient_1989} to highlight the need for high-order methods in
properly capturing the interaction of shocks with smooth wave-like
features. Considering the gathered convergence data, we cannot claim
that the method is of high order away from discontinuities once such
areas enter the domain of influence of a location where artificial
viscosity was applied.  Nonetheless, it is still instructive to see
that the method is capable of keeping the computation stable and
delivering a correct result at least in the ``picture norm'', as
evidenced by Figure \vref{fig:viscous-shu-osher-n5-k80}.  This example
is commonly considered challenging, and it is encouraging that the
method is able to stabilize the computation and give a meaningful
result without excessive smearing.

As a final validation of the detector's design on the Euler equations,
it is important to examine whether it will recognize smooth solutions
and leave them untouched, preserving high-order accuracy. We have
tested this using the smooth isentropic vortex test case of
\citet{zhou_high-resolution_2003} with the result that as soon as
sufficient resolution is available, the detector does not activate
anywhere at any time during the solution process.

\section{Conclusions and Future Work} \label{sec:viscous-conclusions}
What sets the shock detection method of this article apart is its
focus on reliable scaling, with a further emphasis on explicit, local,
GPU-suited calculation in the context of discontinuous Galerkin
methods. Despite a focus on remaining issues, we contend that in this
niche the method is reasonably successful.  Its construction
introduces several new concepts, such as a more precise interpretation
of the correspondence between polynomial decay and smoothness, as well
as methods like skyline pessimization, baseline decay, and $P^1$
viscosity smoothing.

The study of the method's behavior on simple problems (such as linear
waves and transport) was--in our opinion--quite revealing, and it
should be investigated in how far other shock capturing methods are
susceptible to the same problems.

On more complicated nonlinear problems, results were, in our
estimation, encouraging. For example, the method manages to stabilize
the computation of the shock-wave-interaction example and other
important benchmarks, without introducing excessive smoothness.
Further investigation, using the rich pool of tests available in the
shock capturing literature \citep{woodward_numerical_1984,
flash_ug_2009, slater_nparc_2009, stone_test_2009} will doubtlessly
give further insight into the method's strengths and weaknesses as
well as help to further improve it. In addition, we have been
exploring the necessities and pitfalls involved in generalizing the
method to multiple dimensions. Initial tests showed promising results,
which we will report in a future article.

\section*{Acknowledgments}

The authors would like to thank Benjamin Stamm and Gregor Gassner for
valuable discussions, as well as Hendrik Riedmann for contributions to
implementation aspects of this work.  We would also like to thank
Nvidia Corporation for generous hardware donations used to carry out
this research.

TW acknowledges the support of AFOSR under grant number
FA9550-05-1-0473 and of the National Science Foundation under grant
number DMS 0810187. JSH was partially supported by AFOSR, NSF, and
DOE.  AK's research was partially funded by AFOSR FA9550-07-1-0422 and
also through the AFOSR/NSSEFF Program Award FA9550-10-1-0180.  The
opinions expressed are the views of the authors.  They do not
necessarily reflect the official position of the funding agencies.

\footnotesize
\setlength{\bibspacing}{0.3ex}
\bibliographystyle{abbrvnat}
\bibliography{main}

\label{pg:atend}
\end{document}